\documentclass[10pt]{article}
\usepackage{latexsym,amssymb,amsmath}

\def\N{\mathbb{N}}
\def\Q{\mathbb{Q}}
\def\Z{\mathbb{Z}}

\def\C{\mathbb{C}}

\def\proof{\par\medskip\noindent{\em Proof. }}
\def\eproof{\hfill{$\Box$}\bigskip}

\def\tec{\hspace{-1.6mm}{\bf. }}
\def\ds{\dots}
\def\sus{\subset}
\def\al{\alpha}
\def\be{\beta}
\def\ga{\gamma}
\def\de{\delta}
\def\la{\lambda}

\newtheorem{thm}{Theorem}[section]
\newtheorem{prop}[thm]{Proposition}
\newtheorem{cor}[thm]{Corollary}
\newtheorem{lem}[thm]{Lemma}
\newtheorem{prob}[thm]{Problem}
\newtheorem{defi}[thm]{Definition}

\begin{document}
\title{What is an answer?\,---\,remarks, results and problems on PIO formulas in combinatorial 
enumeration, part I}
\author{Martin Klazar}

\maketitle

\begin{abstract}
For enumerative problems, i.e. computable functions $f:\;\N\to\Z$, we define the notion of an {\em 
effective (or closed) formula}.
It is an algorithm computing $f(n)$ in the number of steps that is polynomial in the combined size of the input $n$ 
and the output $f(n)$, both written in binary notation. We discuss many examples of enumerative 
problems for which such closed formulas are, or are not, known. These problems include (i) linear recurrence sequences 
and holonomic sequences, (ii) integer partitions, (iii) pattern-avoiding permutations, (iv) triangle-free graphs and (v) regular graphs. In part I we discuss problems (i) and (ii) and defer (iii)--(v) to part II. Besides other results, 
we prove here that every linear recurrence sequence of integers has an effective formula in our sense. 
\end{abstract}

\section{Introduction}
%sumace pres ciselne rozklady, jake fce davaji PIO (jako treba identicka f(lambda)=1), je to W.-ova ot. 
%co ty trivialni rekurence pro p(n), dve jsou netrivialni. Co treba f(lambda)=#casti lambda^3 a podobne. 
%Treba (-1)^#casti - motivaci je fmle pro # neizom. grafu. Podobne pro mnoz. rozklady. 

We define what is an effective (closed, explicit) formula (solution, algorithm) for a problem in 
enumerative combinatorics or number theory. An {\em enumerative problem} or a {\em counting function} is a computable 
function 
$$
f:\;\N\to\Z
$$ 
(in part II also a computable function $f:\;\{0,1\}^*\to\Z$) given usually via an 
algorithm computing it. The algorithm, often inefficient, usually follows straightforwardly from the statement of 
the problem. Our notation: $|X|$ and $\#X$ denote cardinality of a set $X$, $\N=\{1,2,\ds\}$, $\N_0=\{0,1,\ds\}$, 
$\Z$ is the ring of integers, $\Q$ and $\C$ are, respectively, the fields of rational and complex numbers  
and $\{0,1\}^*$ is the set of finite binary words. The asymptotic symbols 
$o(\cdot)$, $O(\cdot)$, $\Omega(\cdot)$, $\Theta (\cdot)$, $\cdot\ll\cdot$, $\cdot\sim\cdot$ and 
$\mathrm{poly}(\cdot)$ have their usual meaning ($\ll$ is synonymous to $O$ and 
$\mathrm{poly}(x)=O((1+|x|)^d)$ for some $d\in\N$). We will discuss many enumerative problems and effective 
formulas.

Here are some examples of enumerative problems $f:\;\N\to\Z$.
\begin{enumerate}
\item The {\em Catalan numbers} $f(n)=c_n=\frac{1}{n}\binom{2n-2}{n-1}$ which count planted trees 
with $n$ vertices, and many other structures.
\item Let $f(n)=c_n$ for even $c_n$ and $f(n)=1$ for odd $c_n$. 
\item A linear recurrence sequence $f(n+k)=\sum_{i=0}^{k-1}a_if(n+i)$, 
given by initial values $f(1),\ds,f(k)\in\Z$ and recurrence coefficients $a_0,\ds,a_{k-1}\in\Z$.
\item We may be interested in $f(n)=\sum_{\lambda\in P(n)}\left(\|\lambda\|^{\|\lambda\|^{\|\lambda\|}}-
\lfloor\log \|\lambda\|\rfloor\right)$ where $\lambda$ runs through all partitions of $n$ and $\|\lambda\|$ 
is the number of parts.
\item Or in the surplus $f(n)$ of the partitions of $n$ into an even number of distinct parts 
from $\{1_1,1_2,2_1,2_2,3_1,\ds\}$ (two-sorted natural numbers) over those with an odd number of distinct parts, 
that is, $f(n)$ is the coefficient of $q^n$ in the expansion of $\prod_{k=1}^{\infty}(1-q^k)^2$. 
\item Another function $f(n)$ counts $1324$-avoiding permutations $a_1a_2\ds a_n$ of $[n]=\{1,2,\ds,n\}$;  
the avoidance means that no four indices $1\le i_1<i_2<i_3<i_4\le n$ exist with $a_{i_1}<a_{i_3}<a_{i_2}<a_{i_4}$.
\item Or $f(n)$ may be the number of labeled triangle-free graphs on the vertex set $[n]$. 
\item In the binary words setting, if $\lambda=0^{m_0}1^{m_1}\ds(n-1)^{m_{n-1}}$, $m_i\in\N_0$, is a 
multiset with $m_0+m_1+\ds+m_{n-1}=n$, then $f(\lambda)\in\N_0$ counts the labeled simple graphs $G$ on the 
vertex set $[n]$ such that for $i=0,1,\ds,n-1$ exactly $m_i$ vertices of $G$ have degree $i$. We encode the 
input $\lambda$ as an element of $\{0,1\}^*$ in an appropriate way which we will discuss later in part II. 
\end{enumerate}
This shows the variety of problems in enumeration one can consider and investigate. We say something 
on each of them from the perspective of Definition~\ref{pio_form}. Here we consider Examples 1--5 and defer the 
remaining ones to part II. 

We put forward our definition of an effective formula for an enumerative problem. The acronym PIO stands 
for {\em polynomial in input and output}.  
\begin{defi}[PIO formula]\tec\label{pio_form}
For a counting function $f:\;\N\to\Z$, a PIO formula is an algorithm, called a PIO algorithm, that for some 
constants $c,d\in\N$ for every input $n\in\N$ computes the output $f(n)\in\Z$ in at most $c\cdot m(n)^d=O(m(n)^d)=
\mathrm{poly}(m(n))$ steps, where
$$
m(n)=m_f(n):=\log(1+n)+\log(2+|f(n)|)
$$
measures the combined complexity of the input and the output. Similarly for counting function $f:\;X\to\Z$ 
defined on a subset $X\sus\N$. 
\end{defi}
We think this is the precise and definitive notion of a ``closed formula'', and the yardstick 
one should use, possibly with some ramifications or weakenings,  to determine if a solution to an enumerative 
problem is effective. We call counting functions possessing PIO formulas shortly {\em PIO functions}.
Definition~\ref{pio_form} repeats (more explicitly) the proposal made already in M. Klazar \cite[p. 10]{klaz10} 
in 2010, the innovation is that meanwhile we learned that the relevant complexity class PIO exists for a long time 
in the literature. In fact, after submitting this text for publication we found out that J. Shallit mentioned 
briefly Definition~\ref{pio_form} as a ``very good formula'' in \cite[slide 3]{shal} in 2016. 

We comment on the definition. The steps mean steps of the formal specification of an algorithm as a multitape 
Turing machine. We are primarily interested in the {\em bit complexity} but sometimes consider 
also the {\em algebraic complexity}, 
the number of required arithmetic operations. We do not consider the {\em space complexity} which concerns 
memory requirements. The shifts $1+n$ and $2+|f(n)|$ serve for removing arguments $0$ and $1$, inconvenient for 
logarithms. The two natural logarithms come 
from the decadic or binary (but not unary!) encoding of numbers: $\log(2+|f(n)|)=\Theta(r)$ where $r$ is 
the number of bits in the binary code for $f(n)\in\Z$. The rationale behind the definition is that any 
algorithm solving a nontrivial enumerative problem $f:\;\N\to\Z$ needs minimum roughly $m(n)$ steps 
just for reading the input $n$ and printing the output $f(n)$, and an effective algorithm takes only 
polynomially many steps in this minimum. We include the output in the complexity of the problem because 
in enumeration typically the output is much larger than 
the input, and thus considering only the input complexity (and ignoring the time it takes to print the output) 
is bound to lead to confusion. But we will see, and a moment of reflection reveals it, that not large outputs 
but on the contrary the unexpectedly small ones 
pose difficulty\,---\,for them one has much less time for effective computation. We reflect such 
``cancellative'' problems by selecting $\Z$, and not $\N$ or $\N_0$, as the codomain of counting functions. 
In Example 1, $m(n)=\Theta(n)$ because $\log(2+c_n)=\Theta(n)$, and in Example 2, $m(n)=\Theta(n)$ for even 
$c_n$ and $m(n)=\Theta(\log(1+n))$ when $c_n$ is odd.

We state the definition of the complexity class PIO (\cite{yana, gure_shel, zoo})
which we above specialized to counting functions; $|u|$ denotes the length $n$ of a binary word 
$u=a_1a_2\ds a_n\in\{0,1\}^*$. 
\begin{defi}[complexity class PIO]\tec\label{pio_func}
A function $f:\;\{0,1\}^*\to\{0,1\}^*$ is in the complexity class PIO if there is 
an algorithm that for every input $u\in\{0,1\}^*$ computes the output $f(u)$ in time polynomial in 
$\max(|u|,|f(u)|)$. 
\end{defi}
PIO belongs to standard complexity classes, see the ``complexity zoo'' \cite{zoo}, but it appears not to be widely 
known. It was introduced implicitly by M. Yanakakis \cite[Theorem 5.1, also pages 86 and 93]{yana} and later explicitly and 
independently by Y. Gurevich and S. Shelah \cite{gure_shel}. Researchers in database
theory measure by it complexity of algorithms, see for example S. Cohen, B. Kimelfeld and Y. Sagiv \cite{cohe_kime_sagi}, 
S. Cohen and Y. Sagiv \cite{cohe_sagi}, Y. Kanza and Y. Sagiv \cite{kanz_sagi} or M. Vardi \cite{vard}. We are not aware of any
mention of the class PIO in enumerative combinatorics where we think it has a natural place.

The title alludes to the pioneering work \cite{wilf} of the late H.\,S. Wilf who was the first to ponder in the light
of computational complexity the question what it precisely means to give an effective, or a nontrivial, 
solution\,---\,an answer\,---\,to a problem in enumeration. The first of two of his definitions says that a {\em nontrivial solution} 
(he actually uses the term {\em effective solution}) is an algorithm that computes $f(n)$ for a counting function 
$f:\;\N\to\N_0$ in the number 
of steps that is $o(\mathrm{List}(n))$ where ``$\mathrm{List}(n)=$ the complexity of producing all of the members of the 
set $S_n$ [where $f(n)=|S_n|$], one at a time, by the speediest known method, and counting them.'' 
(\cite[Definition 1 and the preceding sentence]{wilf}). As it depends on the complexity of the current ``speediest 
known method'', it depends on time 
and progress of knowlege, and so it is not really a mathematical definition but more a heuristic 
to measure effectivity of algorithms. In part II we give another specification of this notion. But we have to 
add here that this is not a bug but a feature of the first definition: ``We will see that a corollary of this 
attitude is that our decision as to what constitutes an answer may be time-dependent: as faster algorithms for listing
the objects become available, a proposed formula for counting the objects will have to be comparably faster to 
evaluate.'' (\cite[p. 289]{wilf}).

The second definition of H.\,S. Wilf says that for superpolynomially growing $f(n)$ (``a problem in the class $\nu\pi$'') 
an {\em effective solution} (he actually uses the term  that a problem is {\em $p$-solved}) is an algorithm that 
computes $f(n)$ in $\mathrm{poly}(n)$ steps 
(\cite[Definition 2]{wilf}). We also have to mention that he does not restrict only to bit complexity 
but allows also other measures of complexity, ``such as multiplication or division of numbers in a certain size 
range, or bit operations, or function evaluations, etc.'' (\cite[p. 290]{wilf}). 

Shortcomings of H.\,S. Wilf's second definition, rectified in Definition~\ref{pio_form}, are that it does not 
take into account complexity of the output and restricts only to functions $f$ with superpolynomial growth.
(In the second definition the function $f(n)$ is not bounded from above, but it appears that tacitly it is 
of at most broadly exponential growth, $\log (1+f(n))\ll n^d$.) Thus Example 2 falls outside his framework, which is kind of unsatisfactory, 
and so we sought better definition. H.\,S. Wilf illustrates his two definitions with the function (we quote from \cite{wilf})
\begin{eqnarray*}
f(n)&=&\sum_{\lambda\in P(n)}\frac{2^{g(\lambda)}}{1^{m_1}\cdot m_1!\cdot 2^{m_2}\cdot m_2!\cdot\ldots\cdot
n^{m_n}\cdot m_n!}\ \mbox{ where}\\
g(\lambda)&=&\frac{1}{2}\bigg(\sum_{i,j=1}^n\mathrm{gcd}(i,j)m_im_j-\sum_{k\ge1}m_{2k-1}\bigg)
\end{eqnarray*}
and $\lambda=1^{m_1}2^{m_2}\ds n^{m_n}$ runs through the partitions of $n$\,---\,$f(n)$ counts the unlabeled 
(i.e., nonisomorphic) graphs on the vertex set $[n]$ (G. P\'olya \cite{poly}). Since $f(n)\sim2^{n(n-1)/2}/n!$, 
we have $\mathrm{List}(n)=\Omega(2^{(1+o(1))n^2/2})$, but the formula based on
the displayed sum over $P(n)$ computes $f(n)$ in $O(p(n)n^d)=O(\exp(c\sqrt{n}))$ steps 
(where $c>0$ is a constant and $p(n)=|P(n)|$ is the number of partitions of $n$), which is asymptotically much 
smaller. Thus we have a nontrivial 
solution. An effective solution is in question because no algorithm is known that would 
compute $f(n)$ in $O(n^d)$ steps. H.\,S. Wilf asks if such algorithm exists, but so far his question remains unanswered. 

H.\,S. Wilf's article \cite{wilf} is discussed by D. Zeilberger \cite{zeil} and it gave rise to the notion of a {\em 
Wilfian formula} (also called a {\em polynomial enumeration scheme}), which is an algorithm working in time 
polynomial in $n$, for enumerative problems $n\mapsto f(n)\in\N_0$ of the type 
$\Omega(n^c)=\log(2+f(n))=O(n^d)$ (for some real constants $0<c<d$). It appears in the works on enumeration 
of Latin squares by D.\,S. Stones \cite{ston} 
or permutations with forbidden patterns by V. Vatter \cite{vatt}, B. Nakamura and D. Zeilberger \cite{naka_zeil}, 
B. Nakamura \cite{naka} and others. Recently \cite{wilf} was invoked by V.\,S. Miller \cite{mill} for counting squares 
in $F_2^{n\times n}$ by a nontrivial formula making exponentially many steps (an effective solution is not known) or 
by M. Kauers and D. Zeilberger \cite{kaue_zeil}.

Perhaps our proposal in Definition~\ref{pio_form} has a certain reinventing-the-wheel quality because 
nowadays, unlike in the times of, say, L. Comtet \cite{comt} or J. Riordan \cite{rior} when computational complexity 
did not exist, it is a common knowledge that an effective solution to a problem means, in the first approximation, 
a polynomial time algorithm (see, for example, B. Edixhoven and J.-M. Couveignes \cite{edix_couv} or V. Becher, P.\,A. Heiber and T.\,A. Slaman \cite{bech_heib_slam}). For this common knowledge we are indebted to A. Cobham \cite{cobh}
and J. Edmonds \cite{edmo} in 1965. But, then, it seems not to be a common knowledge that one should consistently include 
the complexity of the output in the complexity of an enumerative problem. Also, one can still read in 
contemporary literature on enumerative combinatorics statements to the effect that there is no precise definition or determination of a 
closed formula or answer to an enumerative problem. For example, M. Aigner \cite[Introduction, p. 1]{aign} writes
that ``There is no straightforward answer as to what ``determining'' a counting function means.'' or 
F. Ardila \cite[Chapter 1 What is a good answer?]{ardi} concludes that ``{\bf So what is a good answer to an 
enumerative problem? }[emphasize in original] Not surprisingly, there is no definitive answer to this question.'' On the other hand, 
the text of P.\,J. Cameron \cite[Chapter 1.3]{came, came_book} contains an interesting discussion of the complexity 
matters which does reflect the output complexity of enumerative problems, but does not result in concrete 
definition of a closed 
formula. We believe that Definition~\ref{pio_form} gives the definitive and more or less straightforward answer.
Recently (I wrote the main bulk of this article in autumn 2016 and add this in March 2018, and now see that it is
even August) the excellent survey \cite{pak} of I. Pak appeared that also tackles the question what is an 
effective formula in enumeration but it changes nothing on our above discussion.

{\bf Content and main results. }In the following two sections we discuss, from the perspective of Definition~\ref{pio_form}, 
in their order the eight examples given at the beginning. At least, this we initially intended but as the text started 
get too long, we decided to split it in two parts. Here we consider Examples 1--5 and defer Examples 6--8 to 
part II. The length of our text is caused only by the great variety of enumerative problems offering  
themselves for investigation. Section 2 deals 
with Examples 1--3. After establishing in Propositions~\ref{ex1} and \ref{cnmod2} PIO formulas for Examples 1 and 2, 
where Example 2 is chosen to illustrate peculiarity of this notion, we prove in Theorem~\ref{LRSisPIO} that every 
integral linear recurrence sequence has 
a PIO formula (but with a non-effective complexity bound). This seems so far not to be reflected in the 
literature. Propositions~\ref{LRSandPS} and \ref{kron}--\ref{dichotomy} gather tools for proving Theorem~\ref{LRSisPIO}. 
We present some problems and results on holonomic sequences which generalize linear recurrence 
sequences. For example, in Problem~\ref{prob_holo} we ask if every holonomic sequence $f:\;\N\to\Z$ is a PIO function. 

In Section 3 we consider enumerative problems inspired by Examples 4 and 5. Propositions~\ref{comb_rec_forpn} and \ref{gen_fun_forpn} revisit efficient
evaluation of the partition function $p(n)$. It still holds from us some secrets, for example, it is not known how to compute 
efficiently the parity of $p(n)$ (Problem~\ref{parity_prob}). In Proposition~\ref{PIOnumparts} we show that 
counting functions 
like Example 4 are PIO functions and in Proposition~\ref{PIOnumdistparts} we prove it for partitions with distinct 
parts. Hence Corollary~\ref{comp_dist_part}: compositions of $n$ with distinct parts are counted by a PIO function. 
Proposition~\ref{easy_PIO_part} is a general result implying that, for example, the partitions of $n$ whose 
multiplicities of parts divide $n$ are counted by a PIO function. If parts are required to divide $n$, PIO formula
seems not to be known. Another corollary is Corollary~\ref{slowgrowth}: if $g:\;\N\to\N$ strictly increases, grows only 
polynomially and is polynomial-time computable then the partitions of $n$ with parts in $\{g(1),\,g(2),\,\ds\}$ are 
counted by a PIO function. Yet another Corollary~\ref{part_dsq} shows that the number of partitions of $n$ into 
distinct squares is a PIO function. Corollary~\ref{rest_mult} gives PIO formulas for counting functions of partitions
with prescribed multiplicities. Problem~\ref{wilfs_prob}, inspired by H.\,S. Wilf \cite{wilf10}, asks if we can 
effectively count
partitions of $n$ with distinct multiplicities of parts. The well known theorem of E.\,T. Bell \cite{bell} (Proposition~\ref{bell_thm}) 
says that partitions of $n$ that take parts from a fixed finite set $A\sus\N$ are counted by a quasipolynomial 
in $n$. In Corollaries~\ref{bell_general} and \ref{bell_general_jeste} we point out that the argument proving 
Proposition~\ref{bell_thm} gives with almost no change more general results.  Corollary~\ref{PIOnumparts_other}
returns to Example 4: if the function $g:\;\N\to\Z$ has a finite support then $f(n)=\sum_{\la\in P(n)}g(\|\la\|)$ 
is a PIO function because it is a quasipolynomial in $n$. We mention further quasipolynomial enumerative results 
on partitions, Proposition~\ref{zeilb_dist_mult} due to D. Zeilberger \cite{zeil12} and Proposition~\ref{andr_beck_qpoly} 
due to G.\,E. Andrews, M. Beck and N. Robbins \cite{andr_beck_robb}. Rather general quasipolynomial enumerative result was 
obtained by T. Bogart, J. Goodricks and K. Woods \cite{boga_good_wood} (Theorem~\ref{bgw_thm}). Problem~\ref{m_ary_probl}
asks if one can effectively count partitions of $n$ into powers of a fixed integer $m\ge2$ and 
Theorem~\ref{thm_pak_yeli} quotes a recent positive resolution of this problem by I. Pak and D. Yeliussizov 
\cite{pak_yeli,pak_yeli1}. Proposition~\ref{penta_nepenta} points out that the two basic cancellative counting 
problems on partitions, $n\mapsto\sum_{\la\in Q(n)}(-1)^{\|\la\|}$ and $n\mapsto\sum_{\la\in P(n)}(-1)^{\|\la\|}$ 
(where $P(n)$ are all partitions of $n$, and $Q(n)$ are those with distinct parts), are both PIO functions. The 
former follows from the pentagonal identity of L. Euler, and the latter from J.\,W.\,L. Glaisher's identity 
\cite{glai1}. In the former 
case we have almost complete cancellation but in the latter case only little cancellation. Problem~\ref{probl_on_cancel} asks 
when such cancellation for $(-1)^{\|\la\|}$-count of partitions occurs. We look at this problem for partitions
into squares and for partitions into parts from $l$-sorted $\N$ (Example 5 is $l=2$), including the case $l=24$ that
gives the Ramanujan tau function\,---\,the last Problem~\ref{probl_on_pmqnl} concerns computation of coefficients in 
powers of R. Dedekind's $\eta$-function.    

When we below state and prove results on PIO functions for enumerative problems, we are not content
with just saying  that a PIO algorithm for the problem exists\,---\,it would be ironic to refer in such a way to 
efficient algorithms\,---\,but we always indicate if and how the PIO algorithm can be constructed from the given data. 
See, for example, Theorem~\ref{LRSisPIO} or Corollary~\ref{slowgrowth}.  
    
%Nezapomenout na muj analogicky vysledek
%pro zobecnena Stirl. cisla 2. druhu.
%Ted odstavec o obecne rekurentni metode pro enumeraci.To asi dat do casti II.

\section{The numbers of Catalan and Fibonacci}

We start with Example 1. A {\em planted tree} is a finite tree with a distinguished vertex, called a root,
and with every set of children of a vertex linearly ordered. Its size is the number of vertices.
(For a long time I used to call this kind of trees {\em rooted plane trees}, which also some literature uses, 
but F. Bergeron, G. Labelle and P. Leroux \cite{berg_labe_lero} showed me that this terminology is imprecise: embedding 
in the plane gives to the children of the root only cyclic, not linear, ordering.)

\begin{prop}\tec\label{ex1}
Let $c_n$ be the $n$-th Catalan number, the number of (unlabeled) planted trees with size $n$. 
Then $n\mapsto c_n$ is a PIO function and the PIO algorithm is given by the recurrence displayed below. 
\end{prop}
\proof
We only need to know the recursive structure of planted trees. There is just one planted tree with size $1$, and 
for $n\ge2$ every planted tree $T$ of size $n$ bijectively decomposes in an ordered pair $(U,V)$ of planted trees with sizes 
adding to $n$; $U$ is the subtree of $T$ rooted in the first child of $T$'s root and $V$ is the rest of $T$. 
Thus the combinatorial recurrence $c_1=1$ and, for $n\ge2$,   
$$
c_n=\sum_{k=1}^{n-1}c_kc_{n-k}\;.
$$
It implies that $c_n\ge2c_{n-1}$ for $n\ge3$ and by induction $c_n\gg 2^n$. On the other hand, by induction 
$c_n\le(n-1)!\le n^n$ 
for every $n\ge1$ and $\log(2+c_n)\ll n^2$ (such crude but easy to obtain bound suffices for our purposes). Hence $n\ll m(n)\ll n^2$ in this enumerative problem and we need to compute $c_n$ in $O(n^d)$ steps. 

We do it on a Turing machine with six tapes $T_1,\ds,T_6$. Recall that elementary school algorithms multiply two
$O(m)$-bit numbers in $O(m^2)$ steps and add them in $O(m)$ steps (see comments below).
Tape $T_1$ stores, in this order, the binary codes for $c_1,c_2,\ds,c_{n-1}$. 
For $k=1,2,\ds,n-1$ we do the following. We find $c_k$ and $c_{n-k}$ on $T_1$ and write them on 
the respective tapes $T_2$ and $T_3$. This costs $O(\log(1+c_1)+\ds+\log(1+c_{n-1}))=O(n^3)$ steps, say. We compute 
in $O(n^4)$ steps the product $c_kc_{n-k}=:s$ and write it on $T_4$. Tape $T_5$ stores $\sum_{i=1}^{k-1}c_ic_{n-i}=:t$. 
In $O(n^2)$ steps we 
compute the sum $s+t$ and store it on $T_6$. We conclude by rewritting in $O(n^2)$ steps the content of $T_5$ with that of 
$T_6$. After the step $k=n-1$, tape $T_5$ contains $t=\sum_{i=1}^{n-1}c_ic_{n-i}=c_n$ and we copy this in $O(n^2)$ steps 
on $T_1$. The computation of $c_n$ from $c_1,c_2,\ds,c_{n-1}$ takes $O(n\cdot n^4)=O(n^5)$ steps. The recurrence, implemented by the six-tape 
Turing machine, computes $c_n$ from beginning in $\sum_{k=2}^nO(k^5)=O(n^6)=O(m(n)^6)$ steps, and $n\mapsto c_n$ is a PIO function.
\eproof

\noindent
We give some comments. If the Turing machine has only one tape, and the two $O(m)$-bit numbers to be added are stored on it 
one after another, it is impossible to add them in $O(m)$ steps as the reading head has to move back and forth between 
them, and one needs $\Theta(m^2)$ steps for addition. This can be proven similarly as F.\,C. Hennie \cite{henn} proved the 
$\Theta(m^2)$ lower bound on recognition of $m$-bit palindromes. Therefore we use multitape Turing machines. The 
cost of adding or multiplying two numbers is not only the cost of the operation but in practice includes also the 
cost of recalling both operands from the memory, and therefore we analyzed above the computation of $c_n$ in more 
details. But these technicalities cause at worst only polynomial slowdown and are not important for our main 
concern that is a purely qualitative alternative: there is a polynomial time PIO algorithm for the enumerative 
problem considered or its existence is not known. 

The Catalan numbers 
$$
(c_n)_{n\ge1}=(1,\,1,\,2,\,5,\,14,\,42,\,132,\,429,\,1430,\,4862,\,\ds)
$$ 
form sequence A000108 in the database OEIS (The Online Encyclopedia of Integer Sequences) \cite{oeis}. 
Using stronger bound $\log(2+c_n)=\Theta(n)$ and faster integer multiplication 
(consult D. Harvey, J. van der Hoeven and G. Lecerf \cite{harv_hoev_lece} and D. Harvey and J. van der Hoeven 
\cite{harv_hoev} for the state of art and history of multiplication algorithms) 
we can evaluate $c_n$ in, say, $O(n^3\log^d(1+n))$ steps. It is not a surprise 
that Catalan numbers can be effectively computed, but it depends on what ``effectively'' precisely means. The combinatorial recurrence computes $c_n$ in $\mathrm{poly}(n)$ arithmetic operations, and to get $\mathrm{poly}(n)$ 
steps we need that all numbers involved in the computation have $O(n^d)$ digits for a fixed $d$. This is ensured by
(i) the bound $c_n\le n^n$ and (ii) the non-negativity of coefficients in the recurrence which entails that the 
result $c_n$ upper bounds every number arising in evaluating the recurrence. 
But we also need that each $c_n$ has $\Omega(n^c)$ digits for a fixed real $c>0$, which is ensured by the bound 
$c_n\gg 2^n$, so that $\mathrm{poly}(n)$ steps means $\mathrm{poly}(m(n))$ steps and we really have an effective 
algorithm. If, say, for infinitely many $n$ we had the bound $c_n=O(1)$ then $\mathrm{poly}(n)$ steps would cease to mean 
an effective algorithm in the sense of Definition~\ref{pio_form} and we would have to try more, as in Example 2. 

To establish qualitatively a PIO formula for $c_n$, the combinatorial 
recurrence suffices and one does not need generating functions or ``advanced'' formulas for $c_n$ like ($n\in\N$)
$$
c_n=\frac{1}{n}\binom{2n-2}{n-1}=(-1)^{n+1}\frac{4^n}{2}\binom{\frac{1}{2}}{n}\ \mbox{ or }\ 
c_{n+1}=\frac{4n-2}{n+1}c_n\;.
$$
The last recurrence gives a more efficient PIO formula than the combinatorial recurrence. 
The Catalan numbers have asymptotics $c_n\sim cn^{-3/2}4^n$ with a constant 
$c>0$. For asymptotic methods in enumeration consult P. Flajolet and R. Sedgewick \cite{flaj_sedg}, and R. Pemantle 
and M.\,C. Wilson \cite{pema_wils} for the multivariate case.

We move to Example 2. As we will see shortly, the case $f(n)=1$ occurs for infinitely many $n$. Then the reader 
probably realizes that even though the two functions $f_c$ and $f_o$, defined as 
$f_c(n):=c_n$ and $f_o(n):=n$ for even $n$ and $f_o(n):=1$ for odd $n$, are PIO functions and compose to the 
counting function $f(n)=f_o(f_c(n))$ of Example 2, composition of their PIO algorithms is not
a PIO algorithm for $f(n)$. It is an algorithm that always does $\Omega(n^d)$ steps, but we need an algorithm that 
for $n$ with odd $c_n$ makes only $O(\log^d(1+n))$ steps. Naturally, we need to determine effectively which numbers 
$c_n$ are odd.
\begin{prop}\tec\label{cnmod2}
Let $c_n$ be the $n$-th Catalan number. Then the function $f:\;\N\to\N$, $f(n)=c_n$ for even $c_n$ and $f(n)=1$
for odd $c_n$, is a PIO function. The PIO algorithm is described below.
\end{prop}
\proof
The combinatorial recurrence for $c_n$ shows that $c_n$ is odd 
iff $n=2^m$ for an $m\in\N_0$: $c_1=1$ is odd, for odd $n>1$ the number 
$c_n=2(c_1c_{n-1}+\ds+c_{(n-1)/2}c_{(n+1)/2})$ is even and, similarly, for even $n$ the number $c_n=c_{n/2}^2+
2(c_1c_{n-1}+\ds+c_{(n-2)/2}c_{(n+2)/2})$ has the same parity as $c_{n/2}$. Thus we effectively compute 
$f(n)$ as follows. For given $n\in\N$ we first in $O(\log^2(1+n))$ steps determine if $n$ is a power of $2$. 
If it is so, we output $1$ in $O(1)$ steps. Else we output, using the PIO formula for $n\mapsto c_n$, in $O(n^6)$ or so steps the value $c_n$. 
This gives a formula for $f(n)$ that for even $c_n$ makes $O(n^6)$ steps and for odd $c_n$ only $O(\log^2(1+n))$ steps,
which is $O(m(n)^6)$ steps for every $n\in\N$.
\eproof

\noindent
Alternatively, we get a PIO algorithm for the $f(n)$ of Proposition~\ref{cnmod2} or, more generally, we determine 
effectively if a fixed $m\in\N$ divides $c_n$ (or, more generally, a hypergeometric term), by means of A.-M. 
Legendre's formula $k=\nu_p(n!)=\sum_{j\ge1}\lfloor n/p^j\rfloor$ for the largest $k\in\N_0$ for which $p^k$ divides 
$1\cdot2\cdot\ldots\cdot n$ (or by means of a generalization of the formula to products of numbers in 
arithmetic progressions). At the close of the section we mention other effective formulas for modular reductions of 
$c_n$ and similar numbers.

Example 2 illustrates the fact that composition of two PIO algorithms need not be a PIO algorithm. In fact, as one expects, 
composition of two PIO functions need not be a PIO function. An example of such functions is easily constructed 
by taking a computable function not in PIO, e.g. a function $f:\;\{0,1\}^*\to\{0,1\}$ in $\mathrm{EXP}\backslash\mathrm{P}$, 
see Ch.\,H. Papadimitriou \cite[Chapter 7.2]{papa}. Y. Gurevich and S. Shelah \cite[Lemma 2.1]{gure_shel} elaborate 
such example.

Example 3 concerns the ubiquitous linear recurrence sequences.
Best known of them are the {\em Fibonacci numbers} $f_n$, given by $f_0=0,f_1=1$, and $f_{n+2}=f_{n+1}+f_n$
for every $n\in\N_0$. The sequence 
$$
(f_n)=(f_n)_{n\ge1}=(1,\,1,\,2,\,3,\,5,\,8,\,13,\,21,\,34,\,55,\,89,\,144,\,\ds)
$$ 
is sequence \cite[A000045]{oeis}. Since $2f_{n-2}\le f_n\le 2f_{n-1}$ for $n\ge3$, we have exponential bounds 
$(\sqrt{2})^n\ll f_n\ll 2^n$, $n\in\N$. The precise asymptotics is $f_n\sim c\phi^n$ where $c>0$ is a constant 
and $\phi=1.61803\ds$ satisfies $\phi^2-\phi-1=0$. 

In general, a {\em linear recurrence sequence in $\Z$}, abbreviated as a {\em LRS}, is a function $f:\;\N\to\Z$ 
determined by $2k$ integers $a_0,\ds,a_{k-1},f(1),\ds,f(k)$ with $k\in\N_0$ and $a_0\ne0$, and the recurrence relation
$$
f(n+k)=a_{k-1}f(n+k-1)+a_{k-2}f(n+k-2)+\ds+a_0f(n),\ n\in\N\;.
$$
(Later we point out that allowing $a_i$ outside $\Z$ does not give new LRS.) Note that we require
$a_0\ne0$ and the recurrence to hold from the beginning. By reverting the recurrence every LRS $f$ extends naturally 
to $f:\;\Z\to\Q$, for example $(1,2,4,8,\ds)$ extends to $f(n)=2^{n-1}$, $n\in\Z$.
Thus $(0,1,1,1,\ds)$ is {\em not} a LRS, as can be seen be reverting the purported recurrence, but 
it is true that this sequence differs from a LRS in just one term. If $k=0$ or if $f(1)=f(2)=\ds=f(k)=0$, we get 
the {\em zero sequence} that has $f(n)=0$ for every $n\in\N$. 
We say that $f$ (more  precisely, the recurrence) has {\em order $k$}. Like the Fibonacci numbers, every LRS has an 
easy exponential upper bound $|f(n)|\le c^n$ for every $n\in\N$, with the constant $c=k\max_i|a_i|\max_{i\le k}|f(i)|$. Lower bound is a different story, see Proposition~\ref{growth_dich}. 
 
The defining recurrence computes every LRS $f(n)$ in $\Theta(n)$ arithmetic operations. We recall, on the Fibonacci 
numbers $f_n$, the well known and beautiful formula (algorithm) based on binary powering 
that computes $f(n)$ in only $\mathrm{poly}(\log n)$ arithmetic operations. By E. Bach and J. Shallit 
\cite[p. 122]{bach_shal} or D.\,E. Knuth \cite[p. 695]{knut}, it appears first in J.\,C.\,P. Miller and D.\,J. Spencer 
Brown \cite{mill_spen}. For $f_n$ the formula reads: 
if $n=\sum_{i=0}^kb_i2^i$ with $b_i\in\{0,1\}$ is the binary expansion of $n\in\N_0$ (where $0=02^0$) then
$$
f_n=(0,1)\cdot
\prod_{i=0}^k\underbrace{\left(\ds\left(\left(\left(\begin{array}{ll}
1&1\\1&0
\end{array}\right)^{b_i}\right)^2\right)^2\ds\right)^2
}_{i\ \mathrm{squarings}}
\cdot\left(\begin{array}{l}
1\\0
\end{array}\right)\;.
$$
Indeed, if $M$ is the stated $2\times 2$ matrix and $F_n$ is the column $(f_{n+1},f_n)^T$ then, since the first 
row of $M$ records the recurrence for $f_n$ and matrix multiplication is associative, we have $MF_n=F_{n+1}$ and
$F_n=M^nF_0$. The power $M^n$ is then computed by repeated squaring from the $b_i$s. Since $k=O(\log(n+1))$, the formula  
computes $f_n$ in only $O(\log(n+1))$ multiplications of two integral $2\times 2$ matrices, so in only $O(\log(n+1))$ arithmetic operations with integers. This easily extends 
to general LRS. For more information on repeated squaring and computing terms in linear recurrence sequences see J. von zur Gathen and J. Gerhard \cite[Chapters 4.3 and 12.3]{gath_gerh}. Recent contribution to the literature on computing 
linear recurrence sequences is S.\,G. Hyun, S. Melczer and C. St-Pierre \cite{hyun_al}.

For qualitative bit complexity such cleverness seems superfluous, the $n$-th term $f(n)$ of a LRS is an $O(n)$ 
digit number (which cannot be printed in fewer steps) and already the defining recurrence computes $f(n)$ in 
$\mathrm{poly}(n)$ steps, which might be viewed as an efficient computation. But not from the point 
of view of Definition~\ref{pio_form}. Since $f(n)$ may have as few as $O(1)$ digits, to get a PIO formula we 
need to locate effectively these small values (cf. Example 2) and compute them in $\mathrm{poly}(\log n)$ steps. 
This is not automatic by the above displayed $\mathrm{poly}(\log n)$ arithmetic operations formula because only 
$\mathrm{poly}(\log n)$ digit numbers may be used, and for general LRS the displayed formula contains negative 
numbers (the result need not upper bound intermediate values, as in 
$1=2^n-(2^n-1)$). It can be done, albeit on the verge of non-effectivity, and we prove the following theorem.

\begin{thm}\tec\label{LRSisPIO}
Every linear recurrence sequence $f:\;\N\to\Z$ in $\Z$ of order $k\in\N_0$ is a PIO function. In the proof we indicate
an algorithm ${\cal A}$ with inputs $(a,b,n)\in\Z^k\times\Z^k\times\N$, where $k\in\N_0$, and outputs in $\Z$ such that 
for every fixed tuple $(a,b)=(a_0,\ds,a_{k-1},f(1),\ds,f(k))$ with $a_0\ne0$, ${\cal A}$ is a PIO algorithm 
computing the LRS $f(n)$ determined by $(a,b)$. 
\end{thm}

\noindent
We will see that the implicit constant in the $O(m(n)^d)$ complexity bound in these PIO algorithms 
is currently non-effective\,---\,at the present state of knowledge we cannot provide for it any specific value. The theorem follows by combining standard results from the theory of linear recurrence sequences, see the monograph \cite{ever_al} 
of G. Everest, A. van der Poorten, I. Shparlinski and T. Ward, but we did not find it mentioned in \cite{ever_al} or anywhere else.

We review tools for the proof of Theorem~\ref{LRSisPIO}. For background on linear recurrence sequences 
see \cite{ever_al}, W.\,M. Schmidt \cite{schm} or R.\,P. Stanley \cite[Chapter 4]{stanEC1}. By $\overline{\Q}\sus\C$ 
we denote the field of {\em algebraic numbers,} consisting of all roots of monic polynomials from $\Q[x]$. 
The subring of {\em algebraic integers} is formed by all roots of monic polynomials from $\Z[x]$. 
A {\em power sum} is an expression
$$
s(x)=\sum_{i=1}^lp_i(x)\alpha_i^x 
$$
where $l\in\N_0$, $\alpha_i\in\overline{\Q}$ are distinct and nonzero numbers, and $p_i\in\overline{\Q}[x]$ 
are nonzero polynomials. The numbers $\alpha_i$ are the {\em roots} of the power sum. 
A sequence {\em $f:\;\N\to\overline{\Q}$ is represented by a power sum $s(x)$} if 
$f(n)=s(n)$ for every $n\in\N$. The {\em empty power sum} with $l=0$ represents the zero sequence.
We consider the more general {\em linear recurrence sequences in $\Q$}, shortly {\em LRS in 
$\Q$}\,---\,they are defined as LRS, only their values and coefficients of defining recurrences lie in the field $\Q$ 
instead of its subring $\Z$. If $f:\;\N\to\Q$ is a LRS in $\Q$ given by a recurrence 
$f(n+k)=\sum_{i=0}^{k-1}a_if(n+i)$, $a_i\in\Q$ and $a_0\ne0$, the {\em recurrence polynomial $p(x)$} is
$$
p(x)=x^k-a_{k-1}x^{k-1}-a_{k-2}x^{k-2}-\ds-a_0\in\Q[x]\;.
$$
If the recurrence for $f$ has minimum order, we denote $p(x)$ by $p_f(x)$ and call it the {\em characteristic polynomial of 
$f$}; in a moment we show that there is always a unique $p_f(x)$. For example, the characteristic polynomial of the zero 
sequence is the constant polynomial $p_{\equiv0}(x)=1$. For a sequence $f:\;\N\to\C$ and a polynomial 
$p(x)=\sum_{i=0}^ka_ix^i\in\C[x]$ we let $pf:\;\N\to\C$ denote the sequence given by 
$pf(n)=\sum_{i=0}^ka_if(n+i)$. We say that {\em $p$ annihilates $f$} if $pf$ is the zero sequence. 
The set of rational polynomials annihilating $f$ is denoted by $V(f)\sus\Q[x]$. Clearly, $V(f)$ consists exactly of all  
rational recurrence polynomials for $f$ (with $a_0=0$ and non-unit leading coefficient allowed). The following results are well known but we prove them here for reader's convenience and as an workout for the author.

\begin{prop}\tec\label{LRSandPS}
Power sums and linear recurrence sequences have the following properties.
\begin{enumerate}
\item Every sequence $f:\;\N\to\overline{\Q}$ has at most one power sum representation. 
\item If $f:\;\N\to\C$ is a sequence then $V(f)$ is an ideal in the ring $\Q[x]$.
\item If $f:\;\N\to\Q$ is a LRS in $\Q$ then $V(f)=\langle p_f(x)\rangle$ for a unique monic polynomial 
$p_f\in\Q[x]$ with $p_f(0)\ne0$. This unique generator $p_f$ of $V(f)$ is called the characteristic polynomial of 
$f$ and gives the unique minimum order rational recurrence for $f$. 
\item A sequence $f:\;\N\to\Q$ is a LRS in $\Q$ if and only if it is represented by a power sum 
$s(x)$. If it is the case, the roots $\alpha_i$ of $s(x)$ are exactly the roots of $p_f(x)$. 
\end{enumerate}
\end{prop}
\proof
1. It suffices to show that no nonempty power sum represents the zero sequence. For
a power sum $s(x)$ as above we define $\deg s(x)=\sum_{i=1}^l(\deg p_i(x)+1)\in\N_0$. Clearly, only the empty power 
sum has degree $0$. For any nonempty $s(x)$ we define the new power sum  
$\Delta s(x)=s(x+1)-\alpha_1s(x)$. Note that, crucially, $\deg\Delta s(x)=\deg s(x)-1$. 
Also, if $s(x)$ represents the zero sequence then so does $\Delta s(x)$, and no $s(x)$ with $\deg s(x)=1$ 
represents the zero sequence. (In fact, if $\deg s(x)=1$ then $s(n)\ne0$ for every $n\in\N$.)
The last three facts together imply that no nonempty power sum represents the zero sequence.

2 and 3. It is easy to see that for any $p,q\in\C[x]$ and any sequence $f:\;\N\to\C$ we 
have $(p+q)f=pf+qf$ and $(pq)f=p(qf)$. Thus $V(f)$ is an ideal in $\Q[x]$. Every ideal in $\Q[x]$ is principal, 
is generated by a single element, because the ring $\Q[x]$ is Euclidean. Requiring the generator monic makes 
it unique because the units of $\Q[x]$ are exactly the nonzero constants. Finally, $p_f(0)\ne0$ because $f$ 
being a LRS in $\Q$ implies that $V(f)$ contains a $p$ with $p(0)\ne0$. 

4. Suppose that $f:\;\N\to\Q$ is a LRS in $\Q$ with minimum order $k$ and characteristic polynomial $p_f(x)$. 
Thus in the ring of formal power series $\Q[[x]]$ we have the equality 
$$
\sum_{n\ge0}f(n)x^n=\frac{q(x)}{q_f(x)}
$$ 
where $q_f(x)=x^kp_f(1/x)$, $q\in\Q[x]$ has degree $<k$ and $q(x)$ and $q_f(x)$ are coprime (by the minimality of $k$). 
The value $f(0)$ is computed from $f(1)$, $f(2),\ds,f(k)$ by the reverted recurrence; now it would be more convenient 
if $f$ had domain $\N_0$ or even $\Z$ but counting functions have domain $\N$. After decomposing the rational 
function $q(x)/q_f(x)$ in partial fractions, expanding them in $\overline{\Q}[[x]]$ in generalized geometric series, 
and comparing coefficients of $x^n$, we get a representation of $f(n)$ by a power sum $s(x)$. The roots of $s(x)$ 
are exactly the roots of $p_f(x)$ because of the coprimality of $q(x)$ and $q_f(x)$.

Let $f:\;\N\to\Q$ be represented by a power sum $s(x)=\sum_{i=1}^lp_i(x)\alpha_i^x$: $f(n)=s(n)$ for every $n\in\N$. 
We may assume that $f$ is not the zero sequence and so $s(x)$ is nonempty. We show that 
$f$ is a LRS in $\Q$. The argument is of interest because of three invocations of Lemma~\ref{on_lin_syst} below and 
because it uses negative $n\in\Z$. If $d=\max_i\deg p_i(x)$ then $s(x)$ is
a $\overline{\Q}$-linear combination of the $t=(d+1)l$ expressions $x^j\alpha_i^x$ for $0\le j\le d$ and $1\le i\le l$. So 
is every shift $s(x+r)$ for $r\in\N$, as can be seen by expanding $(x+r)^j=\sum_{b=0}^j\binom{j}{b}r^{j-b}x^b$ and
$\al_i^{x+r}=\al_i^r\al_i^x$. By Lemma~\ref{on_lin_syst} there exist coefficients $\be_0,\be_1,\ds,
\be_t\in\overline{\Q}$, not all zero, such that
$$
\be_0s(x)+\be_1s(x+1)+\ds+\be_ts(x+t)=0
$$
identically. But we need coefficients not only in $\overline{\Q}$ but in $\Q$. We set 
$$
V=\{(f(n),f(n+1),\ds,f(n+t))\;|\;n\in\N\}\sus\Q^{t+1}\sus\overline{\Q}^{t+1}
$$ and select a maximum subset $B\sus V$ of linearly independent (over $\overline{\Q}$) vectors. 
By Lemma~\ref{on_lin_syst}, $|B|\le t+1$. Every vector $z\in V$ is a 
$\overline{\Q}$-linear combination of the vectors in $B$. If $|B|=t+1$, the matrix whose rows are the linearly independent
vectors $z\in B$ is a square matrix and thus has linearly independent columns. But the system 
$$
z\cdot(x_0,x_1,\ds,x_t)=0,\ z\in B\;, 
$$
has a nontrivial solution $x_i=\be_i\in\overline{\Q}$ which means that the columns are linearly dependent. Hence 
$|B|\le t$. But $B\sus\Q^{1+t}$ and thus by Lemma~\ref{on_lin_syst} this system has a nontrivial solution 
$x_i=\ga_i\in\Q$, with not all $\ga_i$ zero. So 
$$
z\cdot(\ga_0,\ga_1,\ds,\ga_t)=0\ \mbox{ and }\ \ga_0s(n)+\ga_1s(n+1)+\ds+\ga_ts(n+t)=0
$$
for every $z\in V$ and every $n\in\N$. By part 1, then $\ga_0s(x)+\ga_1s(x+1)+\ds+\ga_ts(x+t)=0$ identically 
(the left side is the empty power sum) and the last displayed equality thus holds for every $n\in\Z$.
Let $u\in\N_0$ and $v\in\N_0$ be the respective minimum and maximum index $r$ with $\ga_r\ne0$, and let 
$w=v-u\in\N_0$.
Then
\begin{eqnarray*}
&&\ga_vf(n+w)+\ga_{v-1}f(n+w-1)+\ds+\ga_uf(n)\\
&&=\ga_ts(n-u+t)+\ga_{t-1}s(n-u+t-1)+\ds+\ga_0s(n-u)=0
\end{eqnarray*}
for every $n\in\N$ (the arguments of $s(\cdot)$ may be negative). After dividing by $\ga_v$ and rearranging we see 
that $f$ is a LRS in $\Q$ of order $w$. The roots of $p_f(x)$ and of $s(x)$ coincide by the previously proved opposite implication and by uniqueness of $s(x)$ proved in part 1. 
\eproof

\noindent
Note that nonzero power sums may vanish for infinitely many $n\in\N$, for example $s(x)=1^x+(-1)(-1)^x$ on $2\N$, 
and the first part is therefore more subtle result than for polynomials. Recalling the zeros of $\sin(\pi x)$ we have
$$
\exp(\pi ix)-\exp(-\pi ix)=\exp(\pi i)^x+(-1)\exp(-\pi i)^x=0\ \mbox{ for every }\ x\in\Z
$$
but this, of course, is not a counterexample to the first part (why?). Also, 
$$
f_n=\frac{1}{\sqrt{5}}\left(\frac{1+\sqrt{5}}{2}\right)^n-\frac{1}{\sqrt{5}}\left(\frac{1-\sqrt{5}}{2}\right)^n
$$
is the familiar power sum representation of the Fibonacci numbers. The next lemma, used several times in the previous 
proof, is a well known result from linear algebra. Its proof is left to the interested reader as an 
exercise.

\begin{lem}\tec\label{on_lin_syst}
Let $m,n\in\N$ with $m<n$. Every linear homogeneous system 
$$
a_{j,1}x_1+a_{j,2}x_2+\ds+a_{j,n}x_n=0_K,\ j=1,2,\ds,m\;, 
$$
with $m$ equations, $n$ unknowns $x_i$, and coefficients $a_{j,i}$ in a field $K$ has a nontrivial solution $x_i\in K$ 
with not all $x_i=0_K$.
\end{lem}

Recall that a {\em root of unity} is a number $\alpha\in\C$ such that $\alpha^k=1$ for some $k\in\N$, i.e. $\alpha$ is 
a root of $x^k-1$. The minimum such $k$ is the {\em order} of $\al$. We say that a power sum $s(x)$ is {\em degenerate} 
if some root $\alpha_i$ or some ratio 
$\alpha_i/\alpha_j$ of two roots is a root of unity different from $1$, else $s(x)$ is {\em non-degenerate}. 
So we allow $1$ as a root in a non-degenerate power sum, and empty power sum is non-degenerate. For any sequence 
$f:\;\N\to X$ and numbers $m\in\N$ and $j\in[m]$, the {\em $m$-section} $f_{j,m}:\;\N\to X$ of $f$ is the subsequence of values 
of $f$ on the residue class $j$ mod $m$:
$$
f_{j,m}(n)=f(j+m(n-1)),\ n\in\N\;.
$$
If $f:\;\N\to\overline{\Q}$ is represented by a power sum $s(x)=\sum_{i=1}^lp_i(x)\alpha_i^x$, then 
the $m$-section $f_{j,m}$ is represented by the power sum
$$
s_{j,m}(x)=\sum_{i=1}^l\alpha_i^{j-m}p_i(j-m+mx)(\alpha_i^m)^x=\sum_{i=1}^rq_i(x)\beta_i^x
$$
where $r\le l$ and $\{\beta_i\;|\;i=1,\ds,r\}\sus\{\alpha_i^m\;|\;i=1,\ds,l\}$\,---\,we collect like terms in the middle
expression so that the numbers $\beta_i$ are distinct and the polynomials $q_i(x)$ nonzero. For example, 
the degenerate power sum $s(x)=2^x+(-2)^x$ has $2$-sections $s_{1,2}(x)=0$ (the empty power sum with $r=0$) 
and $s_{2,2}(x)=2\cdot4^x$. It is not hard to prove that for any $m\in\N$, $f:\;\N\to\Q$ is a 
LRS in $\Q$ if and only if every $m$-section $f_{j,m}$ is a LRS in $\Q$.

\begin{prop}\tec\label{kron}
The following holds for roots of unity and power sums.
\begin{enumerate}
\item If $p\in\Z[x]$ is a monic polynomial with $p(0)\ne0$ and every root of $p$ has modulus at most $1$, then 
every root of $p$ is a root of unity.
\item If $s(x)=\sum_{i=1}^lp_i(x)\alpha_i^x$ is a power sum such that every $\alpha_i$ is an algebraic integer,
$|\alpha_i|\le1$ for every $i$, and $s(n)\in\Q$ for every $n\in\N$, then every $\alpha_i$ is a root of unity.
\end{enumerate}
\end{prop}
\proof
1. This is called Kronecker's theorem. See U. Zannier \cite[Theorem 3.8 and Remark 3.10 (i)]{zann} or E. Bombieri and 
W. Gubler \cite[Theorem 1.5.9]{bomb_gubl} or V.\,V. Prasolov \cite[Theorem 4.5.4]{pras}. 

2. By the assumption and part 4 of Proposition~\ref{LRSandPS}, the sequence $f(n)=s(n)$, $n\in\N$, is a LRS in $\Q$. 
By parts 3 and 4 of Proposition~\ref{LRSandPS}, the numbers $\alpha_i$ are exactly the roots of 
the characteristic polynomial $p_f(x)\in\Q[x]$. Since all $\alpha_i$ are algebraic integers, so are the coefficients 
of $p_f(x)$ (by expressing them in terms of the $\alpha_i$s). But this implies that $p_f(x)\in\Z[x]$. 
Using part 1 of the present proposition, we get that all $\alpha_i$ are roots of unity.
\eproof

\noindent 
On the Internet or even in paper literature one can encounter the erroneous claim that if 
$\al\in\overline{\Q}$ with $|\al|=1$ then $\al$ is a root of unity. The number $\frac{4+3i}{5}$ is a counterexample. Part 1
of Proposition~\ref{kron} shows when arguments of this sort are correct. Concerning exponential lower bounds on growth 
of linear recurrence sequences, there is the next deep result.
\begin{prop}\tec\label{growth_dich}
If $f:\;\N\to\overline{\Q}$ is represented by a non-degenerate power sum whose roots have maximum modulus $\beta>1$, then 
for every $\varepsilon>0$ there is an $n_0\in\N$ such that
$$
|f(n)|>\beta^{(1-\varepsilon)n}\mbox{ for every } n>n_0\;.
$$
\end{prop}
\proof
This is \cite[Theorem 2.3]{ever_al} where the proof is omitted. At \cite[p. 32]{ever_al} the result is attributed 
to J.-H. Evertse \cite{ever} and independently A. van der Poorten and H.\,P. Schlickewei \cite{vdpo_schl}. J.-H. Evertse 
\cite[p. 229]{ever} attributes it to A. van der Poorten \cite{vdpo}. See also A. van der Poorten \cite{vdpo84}.
%projit to, kdo to vlastne dokazal
\eproof

We deduce the following growth dichotomy for LRS that effectively separates small and 
large values.

\begin{prop}\tec\label{dichotomy}
Suppose $f:\;\N\to\Z$ is a LRS of order $k\in\N_0$, represented by a power sum $s(x)$. We let $m\in\N$ be the least 
common multiple of the orders of the roots of unity among the roots $\al_i$ of $s(x)$ and their ratios $\al_i/\al_j$, 
and let $J\sus[m]$ be the set of $j\in\N$ for which the power sum $s_{j,m}(x)$ is empty or has the single root $1$.
Then there exist a real constant $c>1$ and an $n_0\in\N$ such that for every $j\in[m]$ the following holds.
\begin{enumerate}
\item If $j\in J$ then $f_{j,m}(n)$ is a rational polynomial in $n\in\N$ with degree less than $k$.
\item If $j\not\in J$ then $|f_{j,m}(n)|>c^n$ for every $n>n_0$.
\end{enumerate}
\end{prop}
\proof
By parts 3 and 4 of Proposition~\ref{LRSandPS}, all roots $\al_i$ of $s(x)$ are algebraic integers. 
Take a $j\in[m]$ and consider the power sum $s_{j,m}(x)=\sum_{i=1}^rq_i(x)\beta_i^x$ representing 
$f_{j,m}(n)$. The $\beta_i$s are $m$-th powers of $\al_i$s and are algebraic integers too. Also, 
$s_{j,m}(x)$ is non-degenerate.  Suppose that $|\beta_i|\le1$ for every $i=1,2,\ds,r$. By part 2 of 
Proposition~\ref{kron} all $\beta_i$ are roots of unity. But then non-degeneracy of $s_{j,m}(x)$ implies that either $r=0$, 
$s_{j,m}(x)$ is empty and $f_{j,m}$ is the zero sequence, or $r=1$, $\beta_1=1$ and $f_{j,m}(n)=q_1(n)\in\Z$ for
every $n\in\N$. It follows that $q_1\in\Q[x]$ and $\deg q_1\le\max_i\deg p_i<k$ ($p_i$ are the polynomials in $s(x)$).
Thus we get the first case with $j\in J$. If $|\beta_i|>1$ for some $i$, Proposition~\ref{growth_dich} gives the second 
case with $j\not\in J$.  
\eproof

\noindent
Unfortunately, currently no proof of Proposition~\ref{growth_dich} is known giving an explicit upper bound on the 
threshold $n_0$, only its existence is proven. Therefore also the $n_0$ of 
Proposition~\ref{dichotomy} is non-effective (we cannot compute it). Effective versions of much weaker inequalities 
are not known. Already T. Skolem \cite{skol}
proved that if $f(n)$ is a nonzero LRS represented by a non-degenerate power sum then $|f(n)|\ge1$ for every $n>n_0$, 
that is, $f(n)=0$ has only finitely many solutions $n\in\N$. To obtain an effective version of this result with 
an explicit upper bound on $n_0$, that is, on the sizes of solutions, is a famous open problem, mentioned for example 
in T. Tao \cite[Chapter 3.9]{tao_str_ran} or in B. Poonen \cite{poon}. Before we turn to the proof of Theorem~\ref{LRSisPIO} 
we state a corollary of Proposition~\ref{dichotomy}. Recall that a sequence
$f:\;\N\to\Z$ is a {\em quasi-polynomial} if for some $m\in\N$ polynomials $q_1,\ds,q_m\in\Q[x]$ we have 
$f_{j,m}(n)=q_j(n)$ for every $j\in[m]$ and $n\in\N$.

\begin{cor}\tec\label{dusl_dich}
If a LRS $f:\;\N\to\Z$ has subexponential growth, $$\limsup_{n\to\infty} |f(n)|^{1/n}\le1\;,$$ then 
$f(n)$ is a quasi-polynomial.
\end{cor}

\noindent
We remark that one can prove Proposition~\ref{dichotomy} and Corollary~\ref{dusl_dich} in  a conceptually simpler 
(but probably not much shorter) way without Kronecker's theorem, using incomensurability of the frequencies of the 
roots of non-degenerate power sums. 
 
\bigskip\noindent
{\bf Proof of Theorem~\ref{LRSisPIO}. }Let $k\in\N_0$ and $2k$ integers $a_0,\ds,a_{k-1}$, $f(1),\ds$, $f(k)$,  
$a_0\ne0$, be given. We describe a PIO algorithm for the LRS $f:\;\N\to\Z$ defined by 
$$
f(n+k)=\sum_{i=0}^{k-1}a_if(n+i)\;.
$$
We take the recurrence polynomial $p(x)=x^k-a_{k-1}x^{k-1}-\ds-a_0$ of $f$, decompose the generating function
$$
\sum_{n\ge0}f(n)x^n=\frac{r(x)}{q(x)}\in\Q(x),\ q(x)=x^kp(1/x)\ \mbox{ and }\ \deg r(x)<k\;,
$$ 
into partial fractions and as in the proof of the first implication in part 4 of Proposition~\ref{LRSandPS}  
determine from them the power sum $s(x)$ representing $f(n)$. From $s(x)$ we determine the number $m$ and set 
$J\sus[m]$ as defined in Proposition~\ref{dichotomy}. For each $j\in J$ we find the polynomial $q_j\in\Q[x]$ such 
that $\deg q_j<k$ and $q_j(n)=f_{j,m}(n)$ for $n=1,2,\ds,k$. This precomputation can be done algorithmicly. 
Now for an input $n\in\N$ we compute the residue $j\in[m]$ of $n$ modulo $m$. If $j\in J$, we output 
$f(n)=q_j((n+m-j)/m)$. If $j\not\in J$, we compute $f(n)$ by the defining recurrence. 

Correctness of the algorithm follows from Proposition~\ref{dichotomy}. We bound its time complexity
in terms of $m(n)$. The precomputation takes $O(1)$ steps and determining $j$ takes $\mathrm{poly}(\log n)$ steps.
If $j\in J$, computing $q_j((n+m-j)/m)=f(n)$ takes $\mathrm{poly}(\log n)$ steps because we do $O(1)$ 
arithmetic operations with $O(\log(1+n))$ digit numbers. If $j\not\in J$, computing $f(n)$ by 
the defining recurrence takes $\mathrm{poly}(n)$ steps because $f(n)$ is an $O(n)$ digit number for every 
$n\in\N$. As for $m(n)$, if $j\in J$ then $f(n)$ is an $O(\log(1+n))$ digit number ($f_{j,m}(n)$ grows only 
polynomially) and $m(n)=\Theta(\log(1+n))$. If $j\not\in J$ then $f(n)$ is an $\Omega(n)$ digit number (by case 2 
of Proposition~\ref{dichotomy}) and $m(n)=\Theta(n)$. No matter if $j\in J$ or not, for every $n\in\N$ the algorithm 
does $\mathrm{poly}(m(n))$ steps and is a PIO algorithm.
\eproof

\noindent
Since the constant in the $\Omega(n)$ lower bound at the end of the proof is non-effective, the complexity bound 
$\mathrm{poly}(m(n))=O(m(n)^d)$ involves a non-effective constant as well.

Before we turn to holonomic sequences, we discuss the effect of domain extension for recurrence coefficients of a LRS. 
In the definition we required them to lie in $\Z$. Could one get more integer-valued sequences if the coefficients 
lie in a larger domain than $\Z$? The answer is no. We already proved in the proof of the second implication in part 4 
of Proposition~\ref{LRSandPS} that if $K\sus L$ is an extension of fields and $f:\;\N\to K$ is a LRS in $L$ then 
$f$ is in fact a LRS in $K$. This folklore result on linear recurence sequences is mentioned for example in M. Stoll \cite[Lemma 3.1]{stol}. Recurrence coefficients outside $\Q$ thus give nothing new. One can also prove that if $f:\;\N\to\Z$ 
is a LRS in $\Q$ (recurrence coefficients lie in $\Q$), then $f$ is in fact a LRS 
(another recurrence exists with coefficients in $\Z$). See R.\,P. Stanley \cite[Problem 4.1 (a)]{stanEC1} for the proof 
by generating functions and references for this result, known as the Fatou lemma. 
 
%We remark that we can extend 
%Theorem~\ref{LRSisPIO} to rational linear recurrence sequences $f:\;\N\to\Q$ (determined by rational recurrence 
%coefficients $a_i\in\Q$); we only replace $2+|f(n)|$ in the definition of $m(n)$ with $1+\max(|p|,|q|)$ 
%where $f(n)=p/q\in\Q$ is in lowest terms, and in the case $\max_i|\beta_i|\le 1$ we use the squaring algorithm (now $f(n)$% %is in general not a rational polynomial, e.g., $f(n)=(1/2)^n$), we leave details to the interested reader. 
%Ne ne, je to slozitejsi. 

One could also try to extend Theorem~\ref{LRSisPIO} to linear recurrence sequences in $\Q$.
For this one extends the codomain of counting functions from $\Z$ to $\Q$ and in the definition of $m(n)$ 
(Definition~\ref{pio_form}) replaces $|f(n)|$ with $\max(|a|,|b|)$ where $f(n)=\frac{a}{b}\in\Q$, $\mathrm{gcd}(a,b)=1$. 
We hope to return to this question later. 

We conclude the section with some results and problems on computing terms in {\em holonomic sequences}. These 
generalize LRS and are also quite common in enumerative combinatorics and number theory. For simplicity we restrict 
to integer-valued sequences. A sequence $f:\;\N\to\Z$ is {\em holonomic} (synonymous 
terms in use are {\em $P$-recursive} and {\em polynomially recursive}) if for some $k$ rational functions 
$a_0,\ds,a_{k-1}\in\Z(x)$, $k\in\N_0$ and $a_0(x)\ne0$, we have 
$$
f(n+k)=a_{k-1}(n)f(n+k-1)+a_{k-2}(n)f(n+k-2)+\ds+a_0(n)f(n)
$$
for every $n>n_0$. Now the recurrence cannot hold in general from the beginning because of possible zeros of the denominators in the $a_i(x)$. Examples of such sequences are $f(n)=n!$ or the Catalan numbers $f(n)=c_n$ (see the 
``advanced'' recurrence for $c_n$). Unfortunately, holonomic sequences lack some analog of the power sum representation; for a form of the matrix exponential representation (used in the matrix formula for the Fibonacci
numbers) see Ch. Reutenauer \cite{reut}. We propose the following problem.

\begin{prob}\tec\label{prob_holo}
Is it true that every holonomic sequence $f:\;\N\to\Z$ is a PIO function?
\end{prob}

\noindent
A. Bostan, P. Gaudry and E. Schost \cite{bost_gaud_scho} give an algorithm computing the $n$-th term of a holonomic sequence
in $O(n^{1/2}\log^d(1+n))$ arithmetic operations.

Since Example 2 and Proposition~\ref{cnmod2} deal with an effective computation of the function 
$n\mapsto c_n\ \mathrm{mod}\ 2$, we mention a problem and some results on effective computation of modular reductions
of holonomic sequences. For a sequence $f:\;\N\to\Z$ and $m\in\N$, the {\em modular reduction} $n\mapsto f(n)\ 
\mathrm{mod}\ m$ has values in the fixed set of residues $[m]$, and so a PIO formula for it means a computation in 
$\mathrm{poly}(\log n)$ steps. Trivially, modular reduction of every LRS is eventually 
periodic and has therefore a PIO formula. As we saw in the proof of Proposition~\ref{cnmod2}, $c_n$ modulo $2$ is not
eventually periodic. We propose the following problem.

\begin{prob}\tec\label{prob_hol_red}
Is it true that for every $m\in\N$ and every holonomic sequence $f:\;\N\to\Z$ its modular reduction $n\mapsto(f(n)\ 
\mathrm{mod}\ m)\in[m]$ is a PIO function, that is, can be computed in $O(\log^d(1+n))$ steps?
\end{prob}

\noindent
The answer is affirmative for {\em algebraic $f$}, that is, if the generating series 
$f(x)=\sum_{n\ge 1}f(n)x^n$
satisfies a polynomial equation, $P(x,f(x))=0$ for a nonzero polynomial $P\in\Z[x,y]$ (it is not hard to
show that every algebraic sequence is holonomic). This applies to the Catalan numbers as $c(x)=\sum_{n\ge1}c_nx^n$
satisfies $c(x)^2-c(x)+x=0$. See A. Bostan, X. Caruso, G. Christol and P. Dumas \cite{bost_caru_chri_duma} 
for fast algorithm 
computing modular reduction of algebraic $f$. In fact, the answer is affirmative for an  even larger subclass of 
holonomic sequences, namely for {\em rational diagonals}. These are sequences $f:\;\N\to\Z$ representable for $n\in\N$ as
$$
f(n)=a_{n,n,\ds,n}\ \mbox{ where }\ 
\sum_{n_1,\ds,n_k\ge1}a_{n_1,\ds,n_k}x_1^{n_1}\ds x_k^{n_k}=\frac{P(x_1,\ds,x_k)}{Q(x_1,\ds,x_k)}
$$ 
for some polynomials $P,Q\in\Z[x_1,\ds,x_k]$, $Q\ne0$ (one can show that every algebraic sequence is a rational diagonal, 
and that every rational diagonal is holonomic). See \cite{bost_caru_chri_duma}, A. Bostan, G. Christol and P. Dumas 
\cite{bost_chri_duma} and E. Rowland and R. Yassawi \cite{rowl_yass} (and some of the references therein) for more information 
on these two results. At the conclusion of \cite{rowl} E. Rowland mentions  that a conjecture of G. Christol \cite{chri} 
implies affirmative answer to Problem~\ref{prob_hol_red} for any at most exponentially growing holonomic sequence. 
See C. Krattenthaler and T.\,W. M\"uller \cite{krat_mull} (and other works of the authors cited therein) for 
another approach to computation of modular reductions of algebraic sequences. 

\section{Integer partitions}

Let us discuss PIO formulas for enumerative problems related to and motivated by the initial Examples 4 and 5. 
A {\em partition $\lambda$} of a number $n\in\N_0$ is a multiset of natural numbers 
summing to $n$. We write partitions in two formats, 
$$
\lambda=1^{m_1}2^{m_2}\ds n^{m_n},\ m_i\in\N_0,\ \mbox{ and }\ \lambda=(\lambda_1\ge\lambda_2\ge\ds\ge\lambda_k),\ 
\lambda_i\in\N,\, k\in\N_0\;. 
$$
So 
$|\lambda|:=n=\sum_im_ii=\sum_i\lambda_i$. 
The numbers $1,2,\ds,n$ and $\lambda_1,\lambda_2,\ds,\lambda_k$ are the {\em parts} of $\lambda$ and the $m_i$ are their 
{\em multiplicities}. We denote the number of parts in $\lambda$ by $\|\lambda\|$, so $\|\lambda\|=m_1+m_2+\ds+m_n=k$. 
The set of all partitions of $n$ is $P(n)$, their number is $p(n)=|P(n)|$, and $P:=\bigcup_{n\ge0}P(n)$.
For the empty partition $\emptyset=()$ we have $|()|=\|()\|=0$, $P(0)=\{()\}$ and $p(0)=1$. Similar
quantities are $Q(n)$, $q(n)$, and $Q$, defined for partitions with distinct parts (all $m_i\le1$, i.e. $\la_i>\la_{i+1}$). 
$P_k(n)$ are the partitions of $n$ with $k$ parts, and similarly $Q_k(n)$ are those with $k$ distinct parts. The sequence  
$p:\;\N\to\N$ is \cite[A000041]{oeis} and begins 
$$
(p(n))_{n\ge1}=(1,\,2,\,3,\,5,\,7,\,11,\,15,\,22,\,30,\,42,\,56,\,77,\,101,\,135,\,176,\,\ds)\;;
$$
$(q(n))_{n\ge 1}=(1,1,2,2,3,4,5,\ds)$ is \cite[A000009]{oeis}. We give two proofs for the 
well known fact that $p(n)$ can be efficiently computed and is a PIO function in our parlance, and so is $q(n)$. 
What would be a non-efficient computation? For example, the ``cave man formula'' in \cite{zeil}: 
$p(n)=\sum_{\lambda\in P(n)}1$. By the multiplicity 
format, $q(n)\le p(n)\le(n+1)^n$. To get a lower bound, for a given $n\in\N$, $n\ge4$, consider the maximum $m\in\N$ with 
$1+2+\ds+m=\binom{m+1}{2}\le\frac{n}{2}-1$. Then $m=\Theta(n^{1/2})$ and
$$
n=\sum_{i\in X}i+\left(n-\sum_{i\in X}i\right),\ X\sus[m]\;,
$$ 
are $2^m$ different partitions in $Q(n)$. So 
$p(n)\ge q(n)\ge 2^m\gg\exp(\Omega(n^{1/2}))$ for $n\in\N$. Thus for the {\em partition function $p(n)$} we have
$n^{1/2}\ll m(n)\ll n^2$ and need to compute $p(n)$ in $\mathrm{poly}(n)$ steps. Asymptotically, 
$$
p(n)\sim\frac{\exp(\pi\sqrt{2n/3})}{4\cdot3^{1/2}\cdot n}\ \mbox{ and }\ 
q(n)\sim\frac{\exp(\pi\sqrt{n/3})}{4\cdot3^{1/4}\cdot n^{3/4}}\ \mbox{ as $n\to\infty$}
$$
(G.\,H. Hardy and S. Ramanujan \cite{hard_rama}, G. Meinardus \cite{mein}, G.\,E. Andrews \cite{andr}, 
V. Kot\v e\v sovec \cite{kote}). Thus, more precisely, $p(n)$ and $q(n)$ have $\Theta(n^{1/2})$ digits and 
$m(n)=\Theta(n^{1/2})$. For the sake of brevity we treat the PIO algorithms and their complexity in  
Propositions~\ref{comb_rec_forpn} and \ref{gen_fun_forpn} below more schematically and do not discuss their
implementation by multitape Turing machines as we did in Proposition~\ref{ex1}. These omitted details could be
easily filled in, and it is easy to see that the deduced polynomiality of algorithms holds true.

\begin{prop}\tec\label{comb_rec_forpn} For $k,n\in\N$ with $1\le k\le n$ let 
$p_k(n)=|P_k(n)|$ be the number of partitions of $n$ with $k$ parts, and  
let $p_k(n)=0$ and $P_k(n)=\emptyset$ if $k>n$. Then for every $n\ge1$ we have $p_n(n)=p_1(n)=1$, and for every 
$n\ge2$ and every $k$ with $1<k<n$ we have
$$
p_k(n)=p_k(n-k)+p_{k-1}(n-1)\;.
$$
Consequently, $p(n)=p_1(n)+p_2(n)+\ds+p_n(n)$ is a PIO function.
\end{prop}
\proof
The values $p_n(n)=p_1(n)=1$ are trivial. 
The displayed recurrence mirrors the set partition $P_k(n)=A\cup B$ where $A$ are the partitions of $n$
with all $k$ parts at least $2$ and $B$ are the remaining partitions with at least one part $1$. Decreasing each part in 
every $\lambda\in A$ by $1$ gives the bijection $A\to P_k(n-k)$ and removing one part $1$ from every $\lambda\in B$
gives the bijection $B\to P_{k-1}(n-1)$, whence the recurrence.

For given input $n\in\N$ we use the recurrence and the initial and border values and generate the array of $O(n^2)$ 
numbers $(p_k(m)\;|\;1\le k\le m\le n)$ in $O(n^2)$ additions. Another $n-1$ additions produce $p(n)$. 
Every number involved in the computation has $O(n^2)$ digits (in fact, $O(n^{1/2})$ digits), 
and therefore the algorithm makes $O(n^4)$ steps (in fact, $O(n^{5/2})$ steps), which is 
$O(m(n)^8)$ steps (in fact, $O(m(n)^5)$ steps) as $m(n)\gg n^{1/2}$. Therefore the stated recurrence schema is a PIO formula for $p(n)$. 
\eproof

\noindent
The interested reader will find in M. Bodirsky, C. Gr\"opl and M. Kang \cite{bodi_grop_kang} a recurrence schema, 
in its priciple similar to the previous one but much more involved in details, that computes in polynomial time the 
number of labeled planar graphs on the vertex set $[n]$; computation of this number for $n=50$ in one hour is reported.

The second proof shows that L. Euler's generating function formula 
$$
\sum_{n\ge0}p(n)q^n=\prod_{k\ge1}\frac{1}{1-q^k}=\prod_{k\ge1}(1+q^k+q^{2k}+\ds)
$$ 
also gives a PIO formula for $p(n)$. Let $[x^n]a(x):=a_n$ if $a(x)=\sum_{n\ge0}a_nx^n$.

\begin{prop}\tec\label{gen_fun_forpn}
Let $m,n\in\N$. The product $ab$ of two polynomials $a,b\in\Z[x]$ such that $\deg a,\deg b\le n$ and each
coefficient in them has at most $m$ digits can be computed in the obvious way in $O(m^2n^3)$ steps. Thus
$$
p(n)=[q^n]\prod_{k=1}^n(1+q^k+q^{2k}+\ds+q^{\lfloor n/k\rfloor k})
$$ 
is a PIO function.
\end{prop}
\proof
Each of the $1+\deg a+\deg b=O(n)$ sums $[x^k]ab=\sum_{i+j=k}[x^i]a\cdot [x^j]b$, $k\le \deg a+\deg b$, has  
$O(n)$ summands, multiplication in each summand takes $O(m^2)$ steps, and each addition costs $O(m+n)$ steps 
(we add two numbers with $\ll m+\log(1+n)\ll m+n$ digits). The list of coefficients of $ab$ is thus computed in 
$$
\ll n(m^2n+(m+n)n)=O(m^2n^3)
$$ 
steps. 

The value $p(n)$ is a coefficient in the product of $n$ polynomials with degrees at most $n$ and coefficients $0$ and $1$. We apply the lemma about the product of two 
polynomials $n-1$ times and each time we multiply two polynomials with degrees at most $n^2$ and with coefficients of 
size $\le p(n^2)$ that have $\ll n^4$ digits. Thus we compute the product of the $n$ polynomials in 
$O(n(n^4)^2(n^2)^3)=O(n^{15})=O(m(n)^{30})$ steps (recall that $m(n)\gg n^{1/2}$) and see that 
$p(n)$ is a PIO function. 
\eproof

Often less elementary recurrences are invoked to efficiently compute  $p(n)$:
$$
p(n)=\sum_{i\ge1}(-1)^{i+1}(p(n-a_i)+p(n-b_i))\ \mbox{ or }\ p(n)=\frac{1}{n}\sum_{i=1}^{n}\sigma(i)p(n-i)
$$
where $n=1,2,\ds$, $a_i=\frac{i(3i-1)}{2}$ and $b_i=\frac{i(3i+1)}{2}$ are so called {\em (generalized) pentagonal 
numbers}, $p(0)=1$, $p(n):=0$ for $n<0$, and $\sigma(n):=\sum_{d\,|\,n}d$ is the sum of divisors function (G.\,E. Andrews \cite{andr}). For more recurrences for $p(n)$ see Y. Choliy, L.\,W. Kolitsch and A.\,V. Sills \cite{chol_koli_sill}. The pentagonal recurrence
yields an algorithm computing the list of values $(p(m)\;|\;1\le m\le n)$ in $O(n^2)$ steps (D.\,E. Knuth 
\cite[Chapter 7.2.1.4, exercise 20]{knut_fasc}) while the algorithm of Proposition~\ref{comb_rec_forpn} 
makes $O(n^{5/2})$ steps. N. Calkin, J. Davis, K. James, E. Perez and C. Swannack \cite[Corollary 3.1]{calk_al} 
give an algorithm producing this list in $O(n^{3/2}\log^2n)$ steps, which is close to optimum complexity (because 
it takes $\Omega(n^{3/2})$ steps just to print it). 

The exponent $30$ in the proof of Proposition~\ref{gen_fun_forpn} is hilarious but the reader understands that we do not 
optimize bounds and instead focus on simplicity of arguments. We can decrease it by computing the product $a(x)b(x)$ more 
quickly but in a less elementary way in $O((mn)^{1+o(1)})$ steps by 
\cite[Chapter 8.4]{gath_gerh}. A clever implementation of 
the Hardy--Ramanujan--Rademacher analytic formula for $p(n)$ (\cite{andr}) by F. Johansson \cite[Theorem 5]{joha} 
computes $p(n)$ in $O(n^{1/2}\log^{4+o(1)}n)=O(m(n)^{1+o(1)})$ steps, again in close to optimum 
complexity. F. Johansson reports computing $p(10^6)$ by his algorithm in milliseconds and $p(10^{19})$ 
in less than 100 hours; see \cite{joha_blog} for his computation of $p(10^{20})$. Proposition~\ref{gen_fun_forpn} 
represents $p(n)$ as a coefficient in a polynomial from $\Z[x]$. 
J.\,H. Bruiner and K. Ono \cite{brui_ono} recently found another (more complicated) representation in this spirit, which has $(1-24n)p(n)$ as the next to leading coefficient in a monic polynomial from $\Q[x]$, see also J.\,H. Bruiner, K. Ono 
and A.\,V. Sutherland \cite{brui_ono_suth}. Computation-wise for $p(n)$ it lags far behind the H--R--R formula 
(\cite{brui_ono_suth}).

There is an extensive literature on modular properties of $p(n)$ (for both meanings of ``modular''), see 
K. Ono \cite{ono,ono_modu} and the references therein. N. Calkin et al. \cite[Theorem 3.1]{calk_al} can generate the list $(p(k)\ \mathrm{mod}\ m\;|\;1\le k\le n)$ for special prime moduli $m$ depending on $n$ 
in $O(n^{1+o(1)})$ steps. But unlike for the Catalan numbers and algebraic sequences, so far we do not know an
efficient way to determine the parity of individual numbers $p(n)$. 

\begin{prob}\tec\label{parity_prob}
Is the parity of $p(n)$, the function $n\mapsto(p(n)\ \mathrm{mod}\ 2)\in[2]$, a PIO function? That is, can one compute 
it in $O(\log^d(1+n))$ steps (bit operations)?
\end{prob}
By \cite{joha_blog}, ``With current technology, the most efficient way to determine $p(n)$ modulo a small integer is to
compute the full value $p(n)$ and then reduce it.'' (cf. the discussion of Example 2). Cannot we do better? The parity 
of $p(n)$ was investigated by T.\,R. Parkin and D. Shanks \cite{park_shan} already in 1967. At the end of their article 
they ask if it can be computed in $O(n)$ steps. 

We generalize Proposition~\ref{comb_rec_forpn} by replacing $1$ in $p(n)=\sum_{\lambda\in P(n)}1$ with a 
positive PIO function of the number of parts. This includes Example 4.

\begin{prop}\tec\label{PIOnumparts}
If $g:\;\N\to\N$ is a PIO function then $f:\;\N\to\N$,
$$
f(n)=\sum_{\lambda\in P(n)}g(\|\lambda\|)\;,
$$
is a PIO function too. Construction of the PIO algorithm for $f$ from that for $g$ is described in the proof.
\end{prop}
\proof
We set $G(n)=\max(g(1),g(2),\ds,g(n))$. Clearly, 
$$
f(n)=\sum_{k=1}^ng(k)p_k(n)\;,
$$
and so
$$
f(n)\ge p(n)+G(n)-1\;.
$$
Thus the combined input and output complexity of $f(n)$ satisfies
$$
m(n)=\log(1+n)+\log(2+f(n))\gg\log(p(n)+G(n))\gg n^{1/2}+\log(2+G(n))\;. 
$$
By the assumption on $g$ we compute the list $(g(k)\;|\;1\le k\le n)$ in $$\ll n(\log(1+n)+\log(2+G(n)))^d$$ steps, 
for a fixed $d\in\N$. By Proposition~\ref{comb_rec_forpn} we compute the list $(p_k(n)\;|\;1\le k\le n)$ in 
$O(n^{5/2})$ steps. The product $g(k)p_k(n)$ is computed by elementary school multiplication in
$$
\ll\log^2p(n)+\log^2(2+G(n))\ll(n^{1/2}+\log(2+G(n)))^2
$$
steps, and this also bounds the cost of each addition in the sum. The displayed sum therefore computes $f(n)$ in 
\begin{eqnarray*}
&&\ll n(\log(1+n)+\log(2+G(n)))^d+n^{5/2}+n(n^{1/2}+\log(2+G(n)))^2\\
&&\ll (n^{1/2}+\log(2+G(n)))^{d+4}\ll m(n)^{d+4}
\end{eqnarray*}
steps. 
\eproof
%Plati to i s oborem hodnot f rovnym N_0 ci Z?
%Neco podobneho s nasobnostmi? 

\noindent
Functions covered by the proposition include $f(n)=p(n)$ for $g(n)=1$ and the contrived counting function $f(n)$ of 
Example 4. For $g(n)=n$ we get the total number of parts in all partitions of $n$, so
$$
f(n)=\sum_{\lambda\in P(n)}\|\lambda\|=\sum_{i=1}^n\tau(i)p(n-i)
$$
is a PIO function. Here $\tau(i)$ denotes the number of divisors of $i$. One can deduce the last sum (which
itself is a PIO formula for $f(n)$, given one for $p(n)$, no matter that we cannot compute 
effectively $i\mapsto\tau(i)$) by differentiating the generating function 
$\sum_{\lambda\in P}y^{\|\lambda\|}x^{|\lambda|}=\frac{1}{(1-yx)(1-yx^2)\ds}$ by $y$ and then setting $y=1$. 

For partitions with distinct parts we have similar results.
\begin{prop}\tec\label{PIOnumdistparts}
If $g:\;\N\to\N$ is a PIO function then $f:\;\N\to\N$,
$$
f(n)=\sum_{\lambda\in Q(n)}g(\|\lambda\|)\;,
$$
is a PIO function too. Construction of the PIO algorithm for $f$ from that for $g$ is described in the proof.
\end{prop}
\proof
Now $f(n)=\sum_{k=1}^ng(k)q_k(n)$ where $q_k(n)=|Q_k(n)|$ is the number of partitions of $n$ with $k$ distinct 
parts. For $q_k(n)$ we have the recurrence schema $q_1(n)=q_n(n)=1$ for every $n\ge1$, $q_k(n)=0$ for $k>n$, and  
$$
q_k(n)=q_k(n-k)+q_{k-1}(n-k)
$$
for $n\ge2$ and $1<k<n$. Compared to Proposition~\ref{comb_rec_forpn}, this differs in the last summand: now the 
set of $\lambda\in Q_k(n)$ with one part $1$ bijectively corresponds to $Q_{k-1}(n-k)$, delete the $1$ and decrease
each of the remaining $k-1$ parts by $1$. We only replace $p(n)$ with $q(n)$ and $p_k(n)$ with $q_k(n)$ and argue 
as in the proof of Proposition~\ref{PIOnumparts}. 
\eproof

\noindent
Now for $g(n)=n$ we get that the total number of parts in all $\lambda\in Q(n)$,
$$
f(n)=\sum_{\lambda\in Q(n)}\|\lambda\|=\sum_{i=1}^n\tau^{\pm}(i)q(n-i)\;,
$$
is a PIO function. Here $\tau^{\pm}(i):=\sum_{d\,|\,i}(-1)^{d+1}$ is the surplus of the odd divisors of $i$ over the 
even ones. The last sum (again by itself a PIO formula for $f(n)$, given one for $q(n)$) follows by the same manipulation with the generating function $\sum_{\lambda\in Q}y^{\|\lambda\|}x^{|\lambda|}=(1+yx)(1+yx^2)\ds$ as before. 
We remark that 
\begin{eqnarray*}
{\textstyle
(\sum_{\lambda\in P(n)}\|\lambda\|)_{n\ge1}}&=&(1,\,3,\,6,\,12,\,20,\,35,\,54,\,86,\,128,\,192,\,\ds)\mbox{ and }\\
{\textstyle(\sum_{\lambda\in Q(n)}\|\lambda\|)_{n\ge1}}&=&(1,\,3,\,3,\,5,\,8,\,10,\,13,\,18,\,25,\,30,\,\ds)
\end{eqnarray*}
are respective sequences
\cite[A006128]{oeis} and \cite[A015723]{oeis}. The first one was investigated by ``Miss S.\,M. Luthra, University of Delhi'' 
\cite{luth} (see p. 485 for the formula with $\tau(n)$), and the second by A. Knopfmacher and N. Robbins 
\cite{knop_robb} (they deduce the formula with $\tau^{\pm}(n)$).

Recall that a {\em composition $c$} of $n\in\N$ is an ``ordered partition'' of $n$, that is, a tuple 
$c=(c_1,c_2,\ds,c_k)\in\N^k$ with $c_1+c_2+\ds+c_k=n$. It is well known and easy to show that there are $2^{n-1}$ 
compositions of $n$. What is the number $f_{cdp}(n)$ of compositions of $n$ with distinct parts? 

\begin{cor}\tec\label{comp_dist_part}
The number $f_{cdp}(n)$ of compositions of $n$ with no part repeated is a PIO function. The PIO algorithm for $f_{cdp}$ 
is described in the proof.
\end{cor}
\proof
The mapping $(c_1,c_2,\ds,c_k)\mapsto(c_{i_1}>c_{i_2}>\ds>c_{i_k})$ sending a composition of $n$ with 
distinct parts to its decreasing reordering is a $k!$-to-$1$ mapping from the set of compositions of $n$
with $k$ distinct parts onto $Q_k(n)$. Thus $f_{cdp}(n)=\sum_{k=1}^nk!q_k(n)$ and the result is an instance of 
Proposition~\ref{PIOnumdistparts} for $g(n)=n!$ (clearly, $n\mapsto n!$ is a PIO function, also see P.\,B. Borwein 
\cite{borw}).
\eproof

\noindent
The sequence $(f_{cdp}(n))_{n\ge1}=(1,1,3,3,5,11,13,19,27,\ds)$ is \cite[A032020]{oeis}. B. Richmond and 
A. Knopfmacher \cite{rich_knop} note that 
$$
f_{cdp}(n)=\exp((1+o(1))(2n)^{1/2}\log n)
$$ and obtain a more precise asymptotics. See the book of S. Heubach and T. Mansour \cite{heub_mans} 
for many more enumeration problems for compositions and words, especially with forbidden patterns.

We pose the following problem.
\begin{prob}\tec
Give general sufficient conditions on functions $g:\;\N\to\Z$ ensuring that 
$$
f(n)=\sum_{\la\in P(n)}g(\|\la\|)\ \mbox{ and }\ f(n)=\sum_{\la\in Q(n)}g(\|\la\|)
$$
are PIO functions. 
\end{prob}

\noindent
Propositions~\ref{PIOnumparts} and \ref{PIOnumdistparts} say that it suffices when $g$ is a positive PIO function, 
but it would be more interesting to have general sufficient conditions allowing negative values of $g$.
Corollary~\ref{PIOnumparts_other} and Proposition~\ref{penta_nepenta} are motivated by this problem too. 

We generalize Proposition~\ref{gen_fun_forpn}. Many enumerative problems on partition, but of course not all, fit in 
the general schema of counting partitions with prescribed parts and multiplicities: for every triple 
$(n,i,j)\in\N^2\times\N_0$ we say if $i^j$, part $i$ with multiplicity $j$, may or may not appear in the 
counted partitions of $n$. If $m(n)\gg n^c$ with $c>0$ for the counting problem, the simple algorithm of 
Proposition~\ref{gen_fun_forpn} gives a PIO formula. We spell it out explicitly. 

\begin{prop}\tec\label{easy_PIO_part}
Suppose that $X\sus\N$ is a set, $g(n,i,j)\in\{0,1\}$ is a function defined for $n,i\in\N$ and $j\in\N_0$ with 
$i,j\le n$ and computable in $\mathrm{poly}(n)$ steps, and that the function $f:\;\N\to\N_0$, defined by
\begin{eqnarray*}
f(n)&=&|\{\lambda=1^{j_1}2^{j_2}\ds n^{j_n}\in P(n)\;|\;\mbox{$g(n,i,j_i)=1$ for every $i\in[n]$}\}|\\
&=&[q^n]\prod_{i=1}^n\sum_{j=0}^ng(n,i,j)q^{ij}\;,
\end{eqnarray*}
grows for $n\in X$ as $f(n)\gg\exp(n^c)$ with a constant $c>0$. Then the restriction $f:\;X\to\N_0$ is a PIO function.
The proof shows that the algorithm for $g$ constructively gives the PIO algorithm for $f$, provided that the 
function $f$ gets as inputs only elements of $X$.
\end{prop}
\proof
In the formula for $f(n)$ we have a product of $n$ polynomials with degrees at most $n^2$ each and with the coefficients 
$0$ and $1$ that can be computed in $\mathrm{poly}(n)$ steps. We argue as in the proof of Proposition~\ref{gen_fun_forpn} and deduce that $f(n)$ can be computed in $\mathrm{poly}(n)$ steps, which means $\mathrm{poly}(m(n))$ steps for $n\in X$ 
by the assumption on growth of $f$.
\eproof

\noindent
This applies to partitions of $n$ into distinct parts, odd parts, squares, etc.
but first we illustrate the proposition with two problems where the condition defining counted partitions 
$\lambda\in P(n)$ depends on $n$\,---\,then our reflex to write a formula for $\sum_{n\ge0}f(n)q^n$, often an infinite product of simple
factors, fails us as it cannot be done. They are the functions $f_m,f_p:\N\to\N$ where $f_m(n)$ (resp. $f_p(n)$) counts partitions of $n$ such that every nonzero multiplicity (resp. every part with nonzero multiplicity) divides $n$. 
Then 
\begin{eqnarray*}
(f_m(n))_{n\ge 1}&=&(1,\,2,\,3,\,5,\,4,\,10,\,6,\,17,\,14,\,26,\,13,\,66,\,19,\,63,\,60,\,\ds)\mbox{ and}\\
(f_p(n))_{n\ge 1}&=&(1,\,2,\,2,\,4,\,2,\,8,\,2,\,10,\,5,\,11,\,2,\,45,\,2,\,14,\,14,\,\ds)
\end{eqnarray*}
are respective sequences \cite[A100932]{oeis} and \cite[A018818]{oeis}. By Proposition~\ref{easy_PIO_part}, 
applied with $X=\N$ and $g(n,i,j)=1$ if $j$ divides $n$ and $g(n,i,j)=0$ else, the function 
$f_m(n)$ is a PIO function (with PIO algorithm given in Proposition~\ref{easy_PIO_part}) because $g(n,i,j)$ 
is computable even in $\mathrm{poly}(\log n)$ steps 
and $f_m(n)\ge q(n)\gg\exp(\Omega(n^{1/2}))$ (by the above lower bound on $q(n)$). The second function $f_p(n)$ was investigated by D. Bowman, P. Erd\"os and A. Odlyzko \cite{bowm_al} who proved that 
$$
(1+o(1))\bigg(\frac{\tau(n)}{2}-1\bigg)\log n<\log f_p(n)<(1+o(1))\frac{\tau(n)}{2}\log n\;,
$$
in fact with stronger bounds in place of the $o(1)$ terms. A PIO formula for $f_p(n)$ is apparently not known 
(in contrast to $f_m(n)$, the growth condition is not satisfied and small values occur). After M. Agrawal, 
N. Kayal and N. Saxena \cite{AKS} we however know
a PIO function $f:\;\N\to\{0,2\}$ with the property that $f(n)=2\iff f_p(n)=2$ for every $n\in\N$\,---\,we can 
efficiently compute $f_p(n)$ for infinitely many $n$, namely the prime numbers.

We turn to the more standard situation when $g(n,i,j)$ does not depend on $n$. Then we easily obtain
Corollary~\ref{slowgrowth} below. For its proof we need the next lemma which is also used in the proof of 
Corollary~\ref{rest_mult}.

\begin{lem}\tec\label{frobe}
Let $a_1,a_2,\ds,a_k\in\N$ be distinct numbers such that if $d\in\N$ divides each $a_i$ then $d=1$.
\begin{enumerate}
\item There is an $n_0$, specified in the proof, such that for every $n\in\N$ with $n>n_0$ the equation 
$$
n=a_1x_1+\ds+a_kx_k
$$ 
has a solution $x_1,\ds,x_k\in\N_0$. 
\item The same holds even if the $k$ numbers $x_i$ are required be distinct.
\end{enumerate}
\end{lem}
\proof
1. The ideal $\langle a_1,\ds,a_k\rangle$ in the ring $\Z$ shows that $1=a_1b_1+\ds+a_kb_k$ for some 
$b_i\in\Z$. It is not hard to see that one may take all $b_i$ with $|b_i|\le A^{k-1}$ if $|a_i|\le A$ for every $i$.
We set $c=a_1\max_i|b_i|$ and $d=\sum_ia_ic$. It follows that $n_0=d-1$ works because any $n\ge d$
is expressed by the nonnegative solution $x_1=c+l+jb_1$, $x_i=c+jb_i$ if $i>1$, for appropriate
$l\in\N_0$ and $j\in\{0,1,\ds,a_1-1\}$.

2. Now we set $c=2a_1\max_i|b_i|$ and $d=\sum_ia_ic(k+1-i)$. Then $n_0=d-1$ works, because any $n\ge d$
is expressed by the nonnegative solution with distinct coordinates $x_1=ck+l+jb_1$, $x_i=c(k+1-i)+jb_i$ if $i>1$, 
again for appropriate $l\in\N_0$ and $j\in\{0,1,\ds,a_1-1\}$.
\eproof

\noindent
The first part of the lemma is well known and is the simplest version of the Frobenius problem, see the book \cite{alfo} 
of J.\,L.\,R. Alfons\'\i n for more information.

We look at {\em restricted partitions} with parts in a prescribed set $A\sus\N$; let $P_A(n)\sus P(n)$ be their 
set and $p_A(n):=|P_A(n)|$. We show that we can count them efficiently if the elements of $A$ can be 
efficiently recognized and $A$ is not too sparse. To be precise, in general we probably only ``can'' count them 
efficiently because the proof relies on quantities $d\in\N$ and $B\sus\N$ that in general appear not to be 
computable.

\begin{cor}\tec\label{slowgrowth}
Suppose that the function $g=g(n):\;\N\to\N$ increases, is computable in $\mathrm{poly}(n)$ steps and grows 
only polynomially, $g(n)<(1+n)^d$ for every $n\in\N$ and a constant $d\in\N$, and define $f(n)=p_A(n)$ for 
$A=\{g(1),g(2),\ds\}$,
$$
\sum_{n\ge0}f(n)q^n=\prod_{i\ge1}\bigg(1+\sum_{j=1}^{\infty}q^{g(i)j}\bigg)=\prod_{i\ge1}\frac{1}{1-q^{g(i)}}\;.
$$
Then $f(n)$ is a PIO function. We show how to compute the PIO algorithm for $f(n)$ from the algorithm 
for $g(n)$ if we are given the number $d=\mathrm{gcd}(A)=\mathrm{gcd}(g(1),g(2),\ds)\in\N$ 
and a finite set $B\sus A$ with $\mathrm{gcd}(B)=d$.
\end{cor}
\proof
Let $g(n)$, $A$, $d$ and $B$ be as stated (it is easy to see from prime factorizations that such finite subset 
$B$ exists) and let
$$
n_0=\max(\{0\}\cup\{n\in d\N\;|\;p_B(n)=0\})\in\N_0\;. 
$$
Then $n_0<\infty$ by part 1 of Lemma~\ref{frobe} applied to the numbers $\frac{1}{d}B$, and in fact we can compute $n_0$ 
from the given $B$. So $f(n)=0$ if $n\not\in d\N$ and for $n\in\N$ we can decide in $\mathrm{poly}(\log n)$ steps 
if $n\in d\N$. Subsets $S\sus\{g(1),g(2),\ds,g(m)\}\backslash B$, where $m\in\N$ is maximum with 
$$
g(1)+g(2)+\ds+g(m)\le n-n_0-d\;,
$$ 
prove that for $n\in d\N$ with $n\ge n_0+d$ one has $f(n)\gg\exp(\Omega(n^{1/(d+1)}))$ (for each $S$ we complete the sum 
of its elements by an appropriate partition with parts in $B$ to a partition of $n$). 

We compute $f(n)$ effectively as follows. For the input $n\in\N$ we check 
in $\mathrm{poly}(\log n)$ steps if $n\not\in d\N$ and if $n\le n_0$. In the former case we output $f(n)=0$ and 
in the latter case we compute $f(n)$ by brute force. If neither case occurs, we have $n\in d\N$ and $n>n_0$ and compute 
$f(n)$ by Proposition~\ref{easy_PIO_part}, applied with $X=d\N\backslash[n_0]$ and $g(n,i,j)$ defined as
$g(n,i,j)=0$ if $j\ge1$ and $i\not\in A$, and $g(n,i,j)=1$ else (it is clear that the assumptions  are satisfied, we can
check if $i\not\in A$ in $\mathrm{poly}(i)$ steps). It follows that this is a PIO algorithm for $f(n)$. Also, we 
have constructed it explicitly from the algorithm for $g(n)$ and the knowledge of $d$ and $B$.
\eproof
% a co kdyz g(n)=2^n atd. 

\noindent
In general the quantities $d$ and $B$ probably are not computable from the algorithm for the function $g(n)$. Hence, probably, the PIO algorithm for $f(n)$ cannot be computed given only the algorithm 
for $g(n)$. For example, we may take $g(n)=n^2$ (so $d=1$ and $B=\{1\}$) and compute the number $f_{sq}(n)$ 
of partitions of $n$ into squares, $\sum_{n\ge0}f_{sq}(n)q^n=\prod_{k\ge1}(1-q^{k^2})^{-1}$. 
By Corollary~\ref{slowgrowth} we get a PIO function
$$
(f_{sq}(n))_{n\ge1}=(1,\,1,\,1,\,2,\,2,\,2,\,2,\,3,\,4,\,4,\,4,\,5,\,6,\,6,\,6,\,8,\,\ds)\;,
$$ 
\cite[A001156]{oeis}. The function $f_{sq}(n)$ was investigated by J. Bohman, C.-E. Fr\"oberg and H. Riesel 
\cite{bohm_frob_ries}. As we showed in the proof, $f_{sq}(n)\gg\exp(\Omega(n^{1/3}))$. More generally, already in 
1934 E.\,M. Wright found in \cite{wrig} the asymptotics for the number $p_{S_k}(n)$ of partitions of $n$ into $k$-th 
powers $S_k=\{n^k\;|\;n\in\N\}$ ($k\in\N$): 
$$
p_{S_k}(n)\sim\frac{\Delta}{(2\pi)^{(k+1)/2}}\cdot\frac{k^{1/2}}{(k+1)^{3/2}}\cdot n^{\frac{1}{k+1}-\frac{3}{2}}\cdot
\exp(\Delta n^{1/(k+1)})
$$
where
$$
\Delta=(k+1)\cdot\big((1/k)\cdot\Gamma(1+1/k)\cdot\zeta(1+1/k)\big)^{1-1/(k+1)}\;.
$$
More recently the asymptotics for $p_{S_k}(n)$ with $k=2$ (i.e., $f_{sq}(n)$) was treated by R.\,C. Vaughan \cite{vaug}, 
for general $k$ by A. Gafni \cite{gafn}, and for $p_A(n)$ with $A$ formed by values of an integral polynomial by 
A. Dunn and N. Robles \cite{dunn_robl}.

Why not partition $n$ into {\em distinct} squares, $f_{dsq}(n):=[q^n]\prod_{k\ge1}(1+q^{k^2})$? The initially 
somewhat dull sequence 
$$
(f_{dsq}(n))_{n\ge1}=(1,\,0,\,0,\,1,\,1,\,0,\,0,\,0,\,1,\,1,\,0,\,0,\,1,\,1,\,0,\,1,\,1,\,0,\,\ds)\;, 
$$
\cite[A033461]{oeis}, eventually takes off (the first $n$ with $f(n)\ge2$ is $n=25$) and R. Sprague \cite{spra} 
proved in 1948 that $n=128$ is the last number with $f_{dsq}(n)=0$. See M.\,D. Hirschhorn \cite{hirs} for ``an almost 
complete proof'' that 
$$
f_{dsq}(n)\sim c_2n^{-5/6}\exp(c_1n^{1/3}) 
$$
where $c_1=3c_3^{2/3}$, $c_2=c_3^{1/3}/\sqrt{6\pi}$, and $c_3=\sqrt{\pi}(2-\sqrt{2})\zeta(3/2)/8$. We have

\begin{cor}\tec\label{part_dsq}
The function $f_{dsq}(n):\;\N\to\N_0$ counting partitions of $n$ into distinct squares, which is the same as 
counting partitions of $n$ such that each part $i\in\N$ has multiplicity either $0$ or $i$, is a PIO function.
The PIO algorithm is described in the proof. 
\end{cor}
\proof
This follows from Proposition~\ref{easy_PIO_part}, applied with $X=\{129,130,\ds\}$ and $g(n,i,j)=1$ iff 
$(i=k^2\,\&\,j=1)$ or $j=0$, if we show that $f_{dsq}(n)\gg\exp(n^c)$ on $X$ for some $c>0$. 

To obtain such lower bound for $f_{dsq}(n)$ we begin with a lemma: for every $j\in[4]$ and $n\ge4$ we have 
$|\{A\sus[n]\;|\;|A|\equiv j\;(\mathrm{mod}\;4)\}|\gg2^n$ because any $B\sus[n-3]$ can be enlarged by adding 
one of $n-2,n-1,n$ to have cardinality $j$ modulo $4$. Now for given $n\in\N$ with 
$n\equiv j\;(\mathrm{mod}\;4)$, $j\in[4]$ and $n>1000$, consider the maximum $m\in\N$ with 
$1^2+3^2+5^2+\ds+(2m-1)^2\le n-4\cdot129$. Clearly, $m=\Theta(n^{1/3})$. For every subset $C\sus[m]$ with 
$|C|\equiv j\;(\mathrm{mod}\;4)$ the sum $S_C=\sum_{i\in C}(2i-1)^2$ is also $j$ modulo $4$. By R. Sprague's theorem
mentioned above and our selection of $m$ we may partition $\frac{n-S_C}{4}\in\N$ into distinct squares, say
$\frac{n-S_C}{4}=x_1^2+\ds+x_k^2$. But then $n=S_C+(2x_1)^2+\ds+(2x_k)^2$ is a partition of $n$ 
into distinct squares, and distinct subsets $C$ yield distinct such partitions. Using the lemma we get that 
$f_{dsq}(n)\gg2^m=\exp(\Omega(n^{1/3}))$ for $n\in X$. 
\eproof

\noindent
As M.\,D. Hirschhorn \cite{hirs} himself admits, also in \cite{hirs_rama}, his proof of the asymptotics for $f_{dsq}(n)$ 
is not complete. 

\begin{prob}\tec\label{prob_on_dsq}
Derive rigorously asymptotics for $f_{dsq}(n)$. 
\end{prob}

\noindent
Could not theta functions tell us something about $f_{sq}(n)$ or $f_{dsq}(n)$? It transpires that $f_{dsq}(n)$ 
(with other partition counting functions) comes up in quantum statistical physis.
M.\,N. Tran, M.\,V.\,N. Murty and R.\,K. Bhaduri \cite{tran_al} and, recently, M.\,V.\,N. Murthy, M. Brack, 
R.\,K. Bhaduri and J. Bartel \cite{murt_al} investigate the asymptotics of $f_{dsq}(n)$. The main term is derived, 
with some further terms in the asymptotic expansion, but it is not clear to us whether these arguments are rigorous, 
and so we leave Problem~\ref{prob_on_dsq} as it is.

We partitioned $n$ into squares with arbitrary multiplicities, why not partition $n$ into arbitrary parts but with square 
multiplicities, $f_{sm}(n):=[q^n]\prod_{i\ge1}(1+\sum_{j\ge1}q^{ij^2})$? In general setting this leads to the next 
counterpart to Corollary~\ref{slowgrowth}. Now we have no restriction on the growth of the set of allowed 
multiplicities, it may very well be finite. 

\begin{cor}\tec\label{rest_mult}
Suppose that $1\le g_1<g_2<\ds$ is a finite or infinite nonempty increasing sequence of natural numbers such that 
$n\mapsto g_n$ is computable in $\mathrm{poly}(n)$ steps and define $f(n)=f_A(n)\in\N_0$ to be the number of partitions 
$\lambda\in P(n)$ with all nonzero multiplicities in $A=\{g_1,g_2,\ds\}$,
$$
\sum_{n\ge0}f(n)q^n=\prod_{i\ge1}\bigg(1+\sum_{j\ge1}q^{ig_j}\bigg)\;.
$$
Then $f(n)$ is a PIO function. We show how to compute the PIO algorithm for $f(n)$ from the algorithm for $g_n$ if we 
are given the number $d=\mathrm{gcd}(A)=\mathrm{gcd}(g(1),g(2),\ds)\in\N$ and a finite set $B\sus A$ with 
$\mathrm{gcd}(B)=d$.
\end{cor}
\proof
Let $g_n$, $A$, $d$ and $B$ be as stated and let
$$
n_0=\max(\{0\}\cup\{n\in d\N\;|\;f_B(n)=0\})\in\N_0\;.
$$ 
Then $n_0<\infty$ by the 
second part of Lemma~\ref{frobe} and we can compute $n_0$ from $B$. So $f(n)=0$ if $n\not\in d\N$, and for 
any $n\in\N$ the membership of $n$ in $d\N$ is decidable in $\mathrm{poly}(\log n)$ steps. We obtain a lower bound for 
$f(n)$ with large $n\in d\N$ as in the proof of Corollary~\ref{part_dsq}. Namely, we may assume that $\frac{g_1}{d}$ 
is odd (there is certainly a $g_k$ with odd $\frac{g_k}{d}$) and take arbitrary $n\in d\N$ with $n>2n_0+1$ 
and $\frac{n}{d}\equiv r$ modulo $2$, $r\in[2]$. Let $m\in\N$ be maximum with $(1+3+5+\ds+(2m-1))g_1\le n-2(n_0+1)$. 
Then $m=\Theta(n^{1/2})$ and there are $\gg 2^m$ subsets $C\sus[m]$ such that $|C|\frac{g_1}{d}\equiv r$ modulo $2$. 
For each $C$ we have $(n-g_1\sum_{i\in C}(2i-1))/2\in d\N\backslash[n_0]$ and this number equals $\sum_{g\in D}i_gg$ for some $D\sus B$ and $|D|$ distinct numbers $i_g\in\N$. But then 
$$
n=\sum_{i\in C}(2i-1)g_1+\sum_{g\in D}(2i_g)g
$$ 
is a partition of $n$ into parts $\{2i-1\;|\;i\in C\}$, each with multiplicity $g_1$, and parts $\{2i_g\;|\;g\in D\}$, 
with respective multiplicities $\{g\;|\;g\in D\}$. All multiplicities are in $A$ and distinct 
subsets $C$ give distinct such partitions. Thus for $n\in d\N$ with $n>n_0$ one has $f(n)\gg2^m=\exp(\Omega(n^{1/2}))$.

We compute $f(n)$ effectively as follows. For the input $n\in\N$ we check in $\mathrm{poly}(\log n)$ steps if  
$n\not\in d\N$ and if $n\le n_0$. In the former case we output $f(n)=0$ and in the latter case we compute $f(n)$ by 
brute force. If neither case occurs, we have $n\in d\N$ and $n>n_0$ and compute $f(n)$ by Proposition~\ref{easy_PIO_part}, applied with $X=d\N\backslash[n_0]$ and $g(n,i,j)$ defined for $j\ge1$ and $j\not\in A$ as $g(n,i,j)=0$ and 
$g(n,i,j)=1$ else (clearly the assumption is satisfied, this $g(n,i,j)$ is computable in $\mathrm{poly}(n)$ steps). 
We have for $f(n)$ a PIO algorithm which we have constructed explicitly from the algorithm for $g_n$ and the knowledge 
of $d$ and $B$.
\eproof

\noindent
Similarly to Corollary~\ref{slowgrowth}, the PIO algorithm for $f(n)$ probably cannot be computed given only the 
algorithm for $g_n$. The sequence
$$
(f_{sm}(n))_{n\ge1}=(1,\,1,\,2,\,3,\,3,\,5,\,6,\,8,\,12,\,12,\,17,\,23,\,27,\,32,\,41,\,52,\,\ds)
$$
corresponding to $g_n=n^2$ ($d=1$ and $B=\{1\}$) was for a long time we were preparing this article absent 
in OEIS, but checking it once more in August 2018 we found out that S. Manyama had added it as 
\cite[A300446]{oeis} in May 2018, thank you.

An interesting counting problem for partitions outside the framework of Proposition~\ref{easy_PIO_part} is the number 
$f_{dm}(n)$ of partitions of $n$ into parts with distinct nonzero multiplicities, so
$$
f_{dm}(n)=[q^n]^*\prod_{i=1}^n\bigg(1+\sum_{j=1}^nq^{ij}\bigg)
$$
where the ``star extraction'' of the coefficient means that in the product we only accept monomials 
$q^{i_1j_1}q^{i_2j_2}\ds q^{i_kj_k}$, $1\le i_1<\ds<i_k\le n$, with $|\{j_1,\ds,j_k\}|=k$. The sequence
$$
(f_{dm}(n))_{n\ge1}=(1,\,2,\,2,\,4,\,5,\,7,\,10,\,13,\,15,\,21,\,28,\,31,\,45,\,55,\,62,\,\ds)
$$ 
is \cite[A098859]{oeis}. The problem to investigate $f_{dm}(n)$ was posed by H.\,S. Wilf \cite{wilf10} in 
2010. To get a lower bound on $f_{dm}(n)$, for given $n\in\N$ consider the maximum $m\in\N$
such that $1m+2(m-1)+\ds+m1=(m+1)\sum_{i=1}^mi-\sum_{i=1}^mi^2\le n$. Then $m=\Theta(n^{1/3})$ and the subsets of 
$\{2^{m-1},3^{m-2},\ds,m^1\}$, where the exponents are used in the multiplicities sense, show that 
$f_{dm}(n)\ge 2^{m-1}\gg\exp(\Omega(n^{1/3}))$ (we can always add at least $m$ $1$s to get a distinct 
multiplicities partition of $n$). Thus $f_{dm}(n)$ still has a broadly exponential growth but it is 
not clear how to compute it efficiently.  

\begin{prob}\tec\label{wilfs_prob}
Is the number $f_{dm}(n)$ of distinct-multiplicity partitions of $n$ a PIO function? That is, can we compute
it in $O(n^d)$ steps for a constant $d\in\N$? 
\end{prob}

\noindent
D. Zeilberger \cite{zeil12} says: ``I conjecture that the fastest algorithm takes exponential 
time, but I have no idea how to prove that claim.'' See \cite[A098859]{oeis} for the first 500 or so terms of 
$(f_{dm}(n))_{n\ge1}$. J.\,A. Fill, S. Janson and M.\,D. Ward \cite{fill_al} proved that 
$f_{dm}(n)=\exp((1+o(1))\frac{1}{3}(6n)^{1/3}\log n)$ and D. Kane and R.\,C. Rhoades \cite{kane_rhoa} obtained an even more 
precise asymptotics.

So far we mostly considered partition counting functions $f(n)$ of broadly exponential growth, satisfying 
$\log(2+f(n))\gg n^c$ for a constant $c>0$. We turn to broadly polynomial growth when $\log(2+f(n))\ll\log^d(1+n)$ 
for a constant $d\in\N$. In Corollary~\ref{slowgrowth} we effectively computed $p_A(n)$ for infinite and not 
too sparse sets $A\sus\N$. At the opposite extreme lie finite sets $A\sus\N$, for them $p_A(n)\ll n^{|A|}$. 
By the classical result of E.\,T. Bell \cite{bell} (going back to J. Sylvester in 1857, as E.\,T. Bell himself acknowledges in 
\cite{bell}), then $p_A(n)$ can also be effectively computed.

\begin{prop}[E.\,T. Bell, 1943]\tec\label{bell_thm}
For any finite set $A\sus\N$, the number $p_A(n)$ of partitions $\la=(\la_1\ge\ds\ge\la_k)\in P(n)$ with 
all $\la_i\in A$ is a rational quasipolynomial in $n$. 
\end{prop}

\noindent
Hence for every finite $A$, $p_A(n)$ is a PIO function. Below we extend Proposition~\ref{bell_thm} and indicate 
how to obtain the PIO algorithm. Recall from the previous section that a quasipolynomial $f:\;\Z\to\C$ 
is determined by a modulus $m\in\N$ and $m$ polynomials $p_i\in\C[x]$, $i\in[m]$, such that $f(n)=p_i(n)$ 
if $n\equiv i$ modulo $m$. If $d\ge\deg p_i(x)$ for every $i\in[m]$, we say that 
the quasipolynomial $f$ {\em has class $(m,d)$.} Replacing in the definition $\C$ with $\Q$, we get {\em rational quasipolynomials}. Equivalently, $f$ is a quasipolynomial if and only if 
$$
f(n)=a_k(n)n^k+\ds+a_1(n)n+a_0(n)
$$ 
for some periodic functions $a_i:\;\Z\to\C,\Q$. If $f(n)=p_i(n)$, $n\equiv i$ 
modulo $m$, only holds for every $n\ge N$ for some $N$, we speak of {\em eventual quasipolynomials}. Particular 
cases are {\em (eventually) periodic} sequences $f:\N_0\to\C,\Q$ which are the constant (eventual) quasipolynomials,
each polynomial $p_i(x)$ is (eventually) constant.

There are many treatments of E.\,T. Bell's result in the literature. We mention only R.\,P. Stanley \cite[Chapter 4.4]{stanEC1}, 
\O.\,J. R\o dseth and J.\,A. Sellers \cite{rods_sell}, M. Cimpoea\c{s} and F. Nicolae \cite{cimp_nico1,cimp_nico2}, the 
inductive proof of R. Jakimczuk \cite{jaki} or the recent S. Robins and Ch. Vignat \cite{robi_vign}. The last reference
gives a very nice, simple and short proof of Proposition~\ref{bell_thm} by generating functions, which moreover
presents a PIO formula for $p_A(n)$, $A=\{a_1,a_2,\ds,a_k\}$, explicitly:
$$
p_A(n)
=\sum_{j\in J,\,a_1j_1+\ds+a_kj_k\equiv n\;(\mathrm{mod}\;D)}\binom{\frac{1}{D}(n-a_1j_1-\ds-a_kj_k)+k-1}{k-1}
$$
where $D$ is any common multiple of $a_1,\ds,a_k$ and $J=[0,\frac{D}{a_1}-1]\times\ds\times[0,\frac{D}{a_k}-1]$. 
We give yet another proof via a closure property. We prove that the family 
of rational quasipolynomials is closed under convolution. The {\em (Cauchy) convolution} of functions 
$f,g:\;\N_0\to\C$ is the function $f*g:\;\N_0\to\C$, 
$$
(f*g)(n)=\sum_{i=0}^nf(i)g(n-i)=\sum_{i+j=n}f(i)g(j)\;. 
$$
Convolution gives coefficients in products of power series, if $a,b\in\C[[x]]$ then 
$[x^k]ab=([x^n]a)*([x^n]b)(k)$. If we can compute a rational quasipolynomial $f(n)$ and know 
its class $(m,d)$, we effectively have a PIO algorithm for $f(n)$. Namely, we compute the first (in fact, any) $d+1$ 
values $f(n)$ in each congruence class $n\equiv i$ modulo $m$ and find by Lagrange interpolation the $m$ polynomials 
$p_i\in\Q[x]$ fitting them. The $p_i(x)$ then constitute the PIO algorithm for $f(n)$. Thus it is useful to know 
how classes of quasipolynomials transform under convolution (and linear combination). For other closure properties
of sequences under the operation of convolution see S.\,A. Abramov, M. Petkov\v sek and H. Zakraj\v sek \cite{abra_al}. 

\begin{prop}\tec\label{qpakonv}
Let $f,g:\;\N_0\to\Q$, $\al,\be\in\Q$ and $N,N'\in\N_0$. 
\begin{enumerate}
\item If $f$ and $g$ are rational quasipolynomials then so is $f*g$. If $f$ and $g$ have classes, respectively, 
$(m,d)$ and $(m',d')$ then $f*g$ has class $(M,d+d'+1)$ where $M$ is any common multiple of $m$ and $m'$.
\item If $f$ and $g$ are eventual rational quasipolynomials then so is $\al f+\be g$. If $f$ and $g$ have classes, respectively, 
$(m,d)$ for $n\ge N$ and $(m',d')$ for $n\ge N'$ then $\al f+\be g$ has class $(M,\max(d,d'))$ for $n\ge\max(N,N')$ 
where $M$ is any common multiple of $m$ and $m'$.
\item If $f$ and $g$ are eventual rational quasipolynomials then so is $f*g$. If $f$ and $g$ have classes, 
respectively, $(m,d)$ for $n\ge N$ and $(m',d')$ for $n\ge N'$ then $f*g$ has class $(M,d+d'+1)$ for $n\ge N+N'$ where 
$M$ is any common multiple of $m$ and $m'$. 
\end{enumerate}
\end{prop}
\proof
1. By \cite[Chapter 4.4]{stanEC1} we have 
$$
\sum_{n\ge0}f(n)q^n=\frac{a(q)}{(1-q^m)^{d+1}}\ \mbox{ and }\ 
\sum_{n\ge0}g(n)q^n=\frac{b(q)}{(1-q^{m'})^{d'+1}}
$$ 
for some polynomials $a,b\in\Q[q]$ with degree smaller than that of the denominator. Using the identity 
$$
1-q^{ke}=(1-q^e)(1+q^e+q^{2e}+\ds+q^{(k-1)e}),\ e,k\in\N\;,
$$
which is the main trick in \cite{robi_vign}, we can write for any common multiple $M$ of $m$ and $m'$ 
the product in the same form
$$
\sum_{n\ge0}f(n)q^n\cdot\sum_{n\ge0}g(n)q^n=\frac{a(q)}{(1-q^m)^{d+1}}\cdot\frac{b(q)}{(1-q^{m'})^{d'+1}}=
\frac{c(q)}{(1-q^M)^{d+d'+2}}
$$
where $c\in\Q[q]$ with $\deg c<M(d+d'+2)$. By expanding the denominator in generalized geometric series
we get the result.

2. This is immediate to see.

3. Let $\sum_{n\ge0}f(n)q^n=a(q)+b(q)$ and $\sum_{n\ge0}g(n)q^n=c(q)+d(q)$ where $a,c\in\Q[[q]]$ have 
sequences of coefficients that are rational quasipolynomials with the stated classes for $n\ge 0$ and 
$b,d\in\Q[q]$ are polynomials with degrees $\deg b<N$ and $\deg d<N'$. Let $M$ be any common multiple of 
$m$ and $m'$. Then $(f*g)(n)$ is the sequence of coefficients in
$$
(a+b)(c+d)=ac+ad+bc+bd\;.
$$
By parts 1 and 2, the sequences of coefficients in the last four summands are eventual rational quasipolynomials 
with classes, respectively, $(M,d+d'+1)$ for $n\ge0$, $(m,d)$ for $n\ge N'$ (the sequence of coefficients 
in $ad$ is a linear combination of $\deg d+1$ shifts, by numbers $<N'$, of that in $a$), $(m',d')$ for $n\ge N$, 
and $(1,0)$ for $n\ge N+N'-1$. The result for $f*g$ follows by applying part 2 thrice. 
\eproof

We generalize Proposition~\ref{bell_thm} as follows.
\begin{cor}\tec\label{bell_general}
Let $k,m\in\N$ and 
$$
g:\;\N\times\N_0\to\{0,1\}
$$ 
be a function such that $g(i,0)=1$ for $i>k$, 
$g(i,j)=0$ for $i>k$ and $j>0$, and for $i\le k$ each of the $k$ $0$-$1$ sequences $(g(i,j))_{j\ge0}$ 
is $m$-periodic. Then the function $f:\;\N\to\N_0$,
\begin{eqnarray*}
f(n)&=&\#\{\la=1^{j_1}2^{j_2}\ds n^{j_n}\in P(n)\;|\;\mbox{$g(i,j_i)=1$ for every $i=1,2,\ds,n$}\}\\
&=&[q^n]\prod_{i=1}^k\sum_{j=0}^{\infty}g(i,j)q^{ij}\;,
\end{eqnarray*}
is a rational quasipolynomial with class $(m\cdot k!, k-1)$. The PIO algorithm for $f(n)$ follows constructively 
from $k,m$ and $g$.

We state the ``eventual'' version. Let $k,m$ and $g$ be as before, with the modification that each
sequence $(g(i,j))_{j\ge0}$, $i\le k$, is $m$-periodic only for $j\ge N$, for some given $N\in\N_0$. 
Then the function $f(n)$, defined as above, is an eventual rational quasipolynomial of class $(m\cdot k!, k-1)$ 
for $n\ge\binom{k+1}{2}N$. The PIO algorithm for $f(n)$ follows constructively from $k,m,N$ and $g$.
\end{cor}
\proof
For each $i\le k$ we set $G_i(q)=\sum_{j=0}^{\infty}g(i,j)q^{ij}$. So $G_i\in\{0,1\}[[q]]$ has 
$im$-periodic sequence of coefficients (i.e., with class $(im,0)$). By $k-1$ applications of part 1 of 
Proposition~\ref{qpakonv},
$$
f(n)=[q^n]\prod_{i=1}^kG_i(q)
$$
is  a rational quasipolynomial with class $(m\cdot k!, k-1)$. In the ``eventual'' version, $G_i\in\{0,1\}[[q]]$ has 
$im$-periodic sequence of coefficients for $n\ge iN$. The result follows by $k-1$ applications of part 3 of 
Proposition~\ref{qpakonv}.
\eproof

\noindent
For finite $A\sus\N$, Proposition~\ref{bell_thm} is the instance with $k=\max(A)$ and $g(i,j)=1$ if and 
only if $i\in A$ or $j=0$. More generally Corollary~\ref{bell_general} implies, for example, eventual 
quasipolynomiality of the numbers of partitions of $n$ with parts in $A=\{3,4,27\}$ and such that $3$ appears an 
even number of times, except that multiplicity $2018$ is not allowed, and the multiplicity of $27$ equals $2$ 
or $7$ modulo $11$. Not much changes in the proof, using again Proposition~\ref{qpakonv}, of the next generalization,
and we leave it as an exercise.

\begin{cor}\tec\label{bell_general_jeste}
Let $k,m,d\in\N$ and 
$$
g:\;\N\times\N_0\to\Z
$$ 
be a function such that $g(i,0)=1$ for $i>k$, $g(i,j)=0$ for $i>k$ and $j>0$, 
and for every $i\le k$ the function $j\mapsto g(i,j)$ is a rational quasipolynomial of class $(m,d)$. Then 
($\la=1^{j_1}2^{j_2}\ds n^{j_n}$)
$$
f(n):=\sum_{\la\in P(n)}\prod_{i=1}^ng(i,j_i)\in\Z
$$
is a rational quasipolynomial with class $(m\cdot k!,k(d+1)-1)$. The PIO algorithm for $f(n)$ follows constructively 
from $k,m,d$ and $g$. 

If each function $g(i,j)$, $i\le k$, 
is an eventual rational quasipolynomial of class $(m,d)$ for $j\ge N$, for some given $N\in\N_0$, then $f(n)$ is 
an eventual rational quasipolynomial with class $(m\cdot k!,k(d+1)-1)$ for $n\ge\binom{k+1}{2}N$. The PIO algorithm for 
$f(n)$ follows constructively from $k,m,d,N$ and $g$.
\end{cor}

\noindent
For example, the weighted number of partitions $\la\in P_A(n)$, $A=\{3,4,27\}$, with the weight of $\la$
equal to $(-1)^mm^5+3m$ where $m$ is the multiplicity of $4$, is a rational quasipolynomial (we leave determination 
of its class as an exercise for the reader). 
%We will treat yet further generalization to several variables elsewhere. 

{\em Support} of a function is the set of arguments where it attains nonzero values. In Propositions~\ref{PIOnumparts} 
and \ref{PIOnumdistparts} we made our life easy by positivity of the function $g(n)$. Another easy case 
occurs when $g(n)$ is almost always zero.
\begin{cor}\tec\label{PIOnumparts_other}
If $g:\;\N\to\Z$ has finite support $S\sus\N$ with $s=\max(S)$ then both functions $f_1,f_2:\;\N\to\Z$,
$$
f_1(n)=\sum_{\lambda\in P(n)}g(\|\lambda\|)\ \mbox{ and }\ f_2(n)=\sum_{\lambda\in Q(n)}g(\|\lambda\|)\;,
$$
are rational quasipolynomials with class $(s!,s-1)$. The PIO algorithm for $f_i(n)$ follows constructively from 
$s$ and $g$. (Recall that $Q(n)$ are the partitions of $n$ with no part repeated.)
\end{cor}
\proof
Using conjugation of partitions, which is the involution
$$
(\la_1\ge\la_2\ge\ds\ge\la_k)\leftrightarrow 1^{\la_1-\la_2}2^{\la_2-\la_3}\ds(k-1)^{\la_{k-1}-\la_k}k^{\la_k}\;,
$$ 
we get the well known identity 
$$
p_k(n)=p_{[k]}(n)-p_{[k-1]}(n)
$$ 
---\,the number of partitions of $n$ with $k$ parts equals the number of those partitions of $n$ 
with parts in $[k]$ that use part $k$. Thus
$$
f_1(n)=\sum_{\lambda\in P(n)}g(\|\lambda\|)=\sum_{k\in S}g(k)(p_{[k]}(n)-p_{[k-1]}(n))\;.
$$ 
By part 1 of Proposition~\ref{qpakonv}, $p_{[k]}(n)$ is a rational quasipolynomial of class $(k!,k-1)$. 
The result follows by part 2 of Proposition~\ref{qpakonv}.

As for $f_2(n)$, we have $f_2(n)=\sum_{k\in S}g(k)q_k(n)$. Conjugation of partitions 
yields the identity 
$$
q_k(n)=p_{[k]}(n)-\bigg|\bigcup_{i=1}^kP_{[k]\backslash\{i\}}(n)\bigg|
$$ 
---\,the number of partitions of $n$ into $k$ distinct parts equals the number of those partitions of $n$ with 
parts in $[k]$ that use each part $1,2,\ds,k$. Applying the principle of inclusion and exclusion, 
we express $f_2(n)$ as an integral linear combination 
of the functions $p_A(n)$ with $A\sus[s]$. The result follows by part 2 of Proposition~\ref{qpakonv}.
\eproof

\noindent
In the literature one can find several interesting enumerative results on partitions involving quasipolynomials. 
We state the results of D. Zeilberger \cite{zeil12} and of G.\,E. Andrews, M. Beck and N. Robbins \cite{andr_beck_robb}, for 
other quasipolynomial results see A.\,D. Christopher and M.\,D. Christober \cite{chri_chri} and V. Jel\'\i nek and 
M. Klazar \cite{jeli_klaz}.

\begin{prop}[D. Zeilberger, 2012]\tec\label{zeilb_dist_mult}
For every finite set $A\sus\N$ the number $f(n)$ of the partitions $\la\in P_A(n)$ that have distinct 
nonzero multiplicities is a quasipolynomial in $n$ (and so a PIO function). 
\end{prop}

\begin{prop}[G.\,E. Andrews, M. Beck and N. Robbins, 2015]\tec\label{andr_beck_qpoly}\\
Let $t\in\N_0$ with $t\ge2$. The number $f(n)$ of $\la=(\la_1\ge\ds\ge\la_k)\in P(n)$ with $\la_1-\la_k=t$ 
is a quasipolynomial in $n$ (and so a PIO function). 
\end{prop}

\noindent
For $t=0$ and $1$ the reader can check that, respectively, $f(n)=\tau(n)$ and $f(n)=n-\tau(n)$ where
$\tau(n)$ is the number of divisors of $n$ (and not the Ramanujan function which we will discuss too).
Sadly, despite their simplicity, we do not know if these are PIO formulas because we do not know how to 
efficiently factorize numbers. In fact, \cite{andr_beck_robb} contains a more general result for prescribed 
differences between parts. 

Quite general result in enumeration and logic involving quasipolynomials was achieved by T. Bogart, J. Goodrick and 
K. Woods \cite{boga_good_wood}. In the statement we extend in the obvious way the notion of eventual quasipolynomial 
to function defined on an eventually periodic subset of $\N_0$.

\begin{thm}[T. Bogart, J. Goodrick and K. Woods, 2017]\tec\label{bgw_thm}
Let $d\in\N$ and 
$$
t\mapsto X_t\sus\Z^d,\ t\in\N\;, 
$$
be a sequence of sets $X_t$ that is defined by a formula 
in $1$-parametric Pressburger arithmetic. Then the set $Y\sus\N$ of $t\in\N$ for which $|X_t|<\infty$ is eventually 
periodic and $f:\;Y\to\N_0$, $f(t)=|X_t|$, is an eventual quasipolynomial (and so a PIO function).
\end{thm}

\noindent
The way of definition of the sequence $t\mapsto X_t$ is that the membership 
$$
(x_1,x_2,\ds,x_d)\in X_t
$$ 
is defined, for some $k\in\N$, by a formula built by logical connectives and quantification of integer variables from 
atomic inequalities of the form 
$$
a_1y_1+a_2y_2+\ds+a_ky_k\le b\ \mbox{ where }\ a_i,b\in\Z[t]
$$ 
(here enters the single parameter $t$ in the problem) and the $y_i$ are integer variables including the $x_i$.
Consult \cite{boga_good_wood} for details and examples and, of course, for the proof. T. Bogart, J. Goodrick, 
D. Nguyen and K. Woods prove in \cite{boga_good_nguy} that for more than one parameter, polynomial-time 
computability disappears (assuming $\mathrm{P}\ne\mathrm{NP}$). We state Propositions~\ref{zeilb_dist_mult} and \ref{andr_beck_qpoly} 
and Theorem~\ref{bgw_thm} in their original
form and so do not indicate how to get PIO algorithms from the given data, but with some effort such extensions 
probably could be obtained from the proofs.  

In Corollary~\ref{bell_general}, for $i\le k$ each $0$-$1$ sequence $(g(i,j))_{j\ge0}$ recording allowed 
multiplicities of part $i$ follows a linear (periodic or eventually periodic) pattern. What happens for, say, 
quadratic patterns? What is the number 
$f_{x^2+2y^2}(n)$ of partitions $\la=1^{x^2}2^{y^2}\in P(n)$, $x,y\in\N_0$, that is, partitions of $n$ into 
parts $1$ and $2$ with square multiplicities? A nice formula exists:
$$
f_{x^2+2y^2}(n)=\frac{\tau_{1,8}(n)+\tau_{3,8}(n)-\tau_{5,8}(n)-\tau_{7,8}(n)+\de}{2}
$$
where $\tau_{i,m}(n)$ counts divisors of $n$ that are $i$ modulo $m$, $\de=1$ if $n$ is a square or twice a square 
and $\de=0$ else. This goes back to P. Dirichlet in 1840, see M.\,D. Hirschhorn \cite{hirs00} for a proof. From the reason we already mentioned we do not know if it is a PIO formula. We do not know how to count efficiently solutions 
of equations like $n=x^2+2y^2$ or $n=x^2+y^2$ ($x,y\in\N_0$ or, more classically, $x,y\in\Z$). But if only one of the 
patterns is quadratic, we can again count efficiently. For example, we can count efficiently solutions of 
$n=x+2y^2$: 
$$
f_{x+2y^2}(n):=\#\{\la=1^{j_1}2^{j_2^2}\in P(n)\;|\;j_i\in\N_0\}=\lfloor\sqrt{n/2}\rfloor+1
$$
is a PIO function (but not a quasipolynomial). We compute it in $\mathrm{poly}(\log n)$ steps as 
follows. To compute integral square root $n\mapsto\lfloor\sqrt{n}\rfloor$ in $\mathrm{poly}(\log n)$ steps, 
initialize $m:=0$, add to $m$ in $m:=m+2^r$ the largest power of two such that $m^2\le n$, and repeat. 
When $m$ cannot be increased by adding a power of two, $m=\lfloor\sqrt{n}\rfloor$. See \cite[Chapter 9.5 
and Exercise 9.43]{gath_gerh} for faster algorithms. We hope to treat generalizations of $f_{x+2y^2}(n)$ elsewhere. 

In the setup of Corollary~\ref{slowgrowth}, we get for $f(n)$ broadly polynomial growth if $g(n)$ grows at 
least exponentially. The classical example is for $m\in\N$ with $m\ge2$ the counting function
$$
f_{mp}(n)=f_{mp}(n,m):=p_{\{1,m,m^2,m^3,\ds\}}(n)
$$ 
counting the partitions of $n$ in powers of $m$, so called {\em $m$-ary partitions}. For $m=2$ 
we get the {\em binary partitions}. Binary partitions with distinct parts are easy to count as
$\prod_{k=0}^{\infty}(1+q^{2^k})=\sum_{n=0}^{\infty}q^n$ (partition theorist's joke). But it appears not easy to 
count effectively general binary partitions or $m$-ary partitions. ``Effectively'' here means, of course, 
in $\mathrm{poly}(\log n)$ steps: $m(n)=\log(1+n)+\log(2+f_{mp}(n))=\Theta(\log^2(1+n))$ because 
K. Mahler \cite{mahl} proved that $f_{mp}(n)=\exp((1+o(1))(1/2\log m)\log^2n)$. More precise asymptotic relations 
were derived by N.\,G. de Bruijn \cite{brui}, C.-E. Fr\"oberg \cite{frob} and others.

\begin{prob}\tec\label{m_ary_probl}
Let $m\in\N$ with $m\ge2$. Is the function 
$${\textstyle
f_{mp}(n)=\#\{(x_i)_{i\ge0}\sus\N_0\;|\;\sum_{i\ge0}x_im^i=n\}
}
$$
counting $m$-ary partitions of $n$ a PIO function? That is, can we compute it in $O(\log^d(1+n))$ 
steps (bit operations) for a fixed $d\in\N$?
\end{prob}

\noindent
An interesting algorithm of V.\,P. Bakoev \cite{bako} suggests that the answer might be positive.

\begin{prop}[V.\,P. Bakoev, 2004]\tec
Let $m\in\N$ with $m\ge2$ be given. There is an algorithm computing $f_{mp}(m^n)$ for every $n\in\N$ in 
$O(n^3)$ arithmetic operations.
\end{prop}

\noindent
From the literature on $m$-ary partitions we further mention T. Edgar \cite{edga}, M.\,D. Hirschhorn and J.\,A. Sellers 
\cite{hirs_sell} and D. Krenn and S. Wagner \cite{kren_wagn} (which deals mostly with $m$-ary compositions). 
It is easy to see that the number $f_{bp}(n):=f_{mp}(n,2)$ of binary partitions of $n$ follows the
recurrence $f_{bp}(0)=1$ and $f_{bp}(n)=f_{bp}(n-1)+f_{bp}(n/2)$ for $n\ge1$ (where $f_{bp}(n/2)=0$ if 
$n/2\not\in\N_0$). The reduction $f_{bp}'(n):=f_{bp}(2n)$ follows the recurrence $f_{bp}'(0)=1$ and ($n\ge1$)
$f_{bp}'(n)=f_{bp}'(n-1)+f_{bp}'(\lfloor n/2\rfloor)$ and forms the sequence \cite[A000123]{oeis},
$$
(f_{bp}(2n))_{n\ge0}=(f_{bp}'(n))_{n\ge0}=(1,\,2,\,4,\,6,\,10,\,14,\,20,\,26,\,36,\,46,\,60,\,74,\,\ds)\;,
$$
investigated by D.\,E. Knuth \cite{knut66} fifty years ago. 

Recently, I. Pak \cite[Theorem 2.19]{pak} has announced 
positive resolution of Problem~\ref{m_ary_probl} in I. Pak and D. Yeliussizov \cite{pak_yeli,pak_yeli1}: 
(we quote verbatim from \cite{pak})

\begin{thm}[I. Pak and D. Yeliussizov]\tec\label{thm_pak_yeli}
Let ${\cal A}=\{a_1,a_2,\ds\}$, and suppose $a_k/a_{k-1}$ is an integer $\ge2$, for all $k>1$. Suppose also that
membership $x\in{\cal A}$ can be decided in $\mathrm{poly}(\log x)$ time. Then $\{p_{{\cal A}}(n)\}$ can be 
computed in time $\mathrm{poly}(\log n)$.
\end{thm}

In conclusion of Section 3 and of our article we turn to cancellative problems related to the initial Example 5. 
Sums of integers with large absolute values still may be small, even $0$. In enumeration it 
means that a formula, effective for nonnegative summands, may no longer be effective (in the sense of 
Definition~\ref{pio_form}) for integral summands, if cancellations occur. In the next proposition we give
both an example and a non-example of such cancellation. The former is a classics but the latter may be not so 
well known. 

\begin{prop}\tec\label{penta_nepenta}
Both functions
$$
q^{\pm}(n):=\sum_{\la\in Q(n)}(-1)^{\|\la\|}\ \mbox{ and }\ p^{\pm}(n):=\sum_{\la\in P(n)}(-1)^{\|\la\|}
$$
are PIO functions. Concretely, 
$$
\sum_{n=0}^{\infty}q^{\pm}(n)q^n=\prod_{k=1}^{\infty}(1-q^k)=1+\sum_{n=1}^{\infty}(-1)^n(q^{n(3n-1)/2}+q^{n(3n+1)/2})
$$
and 
$$
\sum_{n=0}^{\infty}p^{\pm}(n)q^n=\prod_{k=1}^{\infty}\frac{1}{1+q^k}=
\prod_{k=1}^{\infty}(1+(-q)^{2k-1})=1+\sum_{n=1}^{\infty}(-1)^nq_o(n)q^n
$$
where $q_o(n):=\#\{\la\in Q(n)\;|\;\la_i\equiv1\ (\mathrm{mod}\ 2)\}$.
\end{prop}
\proof
The first identity is the famous {\em pentagonal identity} of L. Euler \cite{eule} (or \cite{andr,hard_wrig}). 
Replacing $q(n)$ with $q^{\pm}(n)$ leads to almost complete cancellation to values just
$0$ and $\pm 1$, and $m(n)=\Theta(\log(1+n))$ for $q^{\pm}(n)$. The algorithms of Propositions~\ref{PIOnumdistparts} 
(recurrence schema) and \ref{gen_fun_forpn} (coefficient extraction from a generating polynomial) still work but 
do $\mathrm{poly}(n)$ steps and are not effective for computing $n\mapsto q^{\pm}(n)$. For a PIO formula 
more efficient algorithm is needed. Fortunately, the pentagonal identity provides it. We easily determine 
in $\mathrm{poly}(\log n)$ steps the existence of a solution $i\in\N$ to the equation $n=\frac{i(3i\pm 1)}{2}$ 
and its parity, simply by computing the integral square root as discussed above in connection with $f_{x+2y^2}(n)$.

The second identity, more precisely the middle equality, follows at once from $\frac{1}{1+q^k}=\frac{1-q^k}{1-q^{2k}}$. 
Now the replacement of $p(n)$ with $p^{\pm}(n)$ leads to almost no cancellation because
$$
|p^{\pm}(n)|=q_o(n)=[q^n]\prod_{k=1}^{\infty}(1+q^{2k-1})\sim\frac{\exp(\pi\sqrt{n/6})}{2^{3/2}\cdot6^{1/4}\cdot n^{3/4}}
$$
(\cite[A000700]{oeis}; V. Kot\v e\v sovec \cite[p. 9]{kote}; G. Meinardus \cite[p. 301]{mein}) remains of broadly 
exponential growth. Thus for $p^{\pm}(n)$ we have $m(n)\gg n^{1/2}$ and both algorithms
for $p(n)$ remain effective for $p^{\pm}(n)$. Hence, more easily than for $q^{\pm}(n)$, $p^{\pm}(n)$ is a 
PIO function.
\eproof

\noindent
We had derived the second identity and then we learned in A. Ciolan \cite{ciol} that it is in fact due to 
J.\,W.\,L. Glaisher \cite{glai1}. The above examples lead us to the following question.

\begin{prob}\tec\label{probl_on_cancel}
Find general sufficient conditions on the functions 
$$
a=a_i:\;\N\to\N_0\ \mbox{ and }\ b=b_{i,j,k}:\;\N\times\N_0\times\N\to\{0,1\}
$$ 
ensuring that 
$$
f(n)=[q^n]\prod_{i=1}^{\infty}\prod_{k=1}^{a_i}\sum_{j=0}^{\infty}b_{i,j,k}(-1)^jq^{ij}
$$
is a PIO function. Find asymptotics of $f(n)$. 
\end{prob}

\noindent
Thus $f(n)$ is the $(-1)^{\|\la\|}$-count of the partitions $\la$ of $n$ into parts
$i\in\N$ coming in $a_i$ sorts $(i,1),(i,2),\ds,(i,a_i)$ such that the  part $(i,k)^j$ may appear if 
and only if $b_{i,j,k}=1$.
For $a_i=1$, $b_{i,0,1}=b_{i,1,1}=1$, and $b_{i,j,k}=0$ else we get $q^{\pm}(n)$, and for $a_i=1$, 
$b_{i,j,1}=1$, and $b_{i,j,k}=0$ else we get $p^{\pm}(n)$. For $a_i=2$, $b_{i,0,1}=b_{i,1,1}=b_{i,0,2}=b_{i,1,2}=1$, 
and $b_{i,j,k}=0$ else we get the function $f(n)$ of Example 5. More generally, for $l\in\N$ the counting functions 
$$
\sum_{n=0}^{\infty}q^{\pm,l}(n)q^n:=\prod_{n=1}^{\infty}(1-q^n)^l
$$
correspond to $a_i=l$, $b_{i,0,1}=b_{i,1,1}=b_{i,0,2}=b_{i,1,2}=\ds=b_{i,0,l}=b_{i,1,l}=1$, and $b_{i,j,k}=0$ else;
$q^{\pm}(n)=q^{\pm,1}(n)$ and Example 5 is $q^{\pm,2}(n)$. A related open problem, due to D. Newman, is mentioned
in G.\,E. Andrews and D. Newman \cite{andr_newm}: In
$$
\sum_{n\ge0}p(n)q^n=\prod_{k=1}^{\infty}(1+p^k+p^{2k}+\ds)\;,
$$
can one change some signs in the last product so that on the left side the $p(n)$ turn to coefficients
$0$ and $\pm1$?

Another example of non-cancellation in Problem~\ref{probl_on_cancel} is the result of A. Ciolan \cite[(22)]{ciol}: 
if $S_2:=\{1,4,9,16,\ds\}$, $B:=\Gamma(3/2)\zeta(3/2)/2\sqrt{2}$ and
$$
t_n:=\sum_{\la\in P(n),\,\la_i\in S_2}(-1)^{\|\la\|}=[q^n]\prod_{k\ge1}\frac{1}{1+q^{k^2}}
$$
then
$$
t_n\sim(-1)^n\frac{\exp\left(3(B/2)^{2/3}n^{1/3}\right)}{(3\pi)^{1/2}(2n)^{5/6}/B^{1/3}}\;.
$$
In \cite[A292520]{oeis}, for example
$$
(t_n)_{n=32}^{49}=(1,\,-2,\,3,\,-4,\,3,\,-2,\,1,\,0,\,1,\,-2,\,3,\,-4,\,3,\,-2,\,1,\,0,\,0,\,-2)\;,
$$
V. Kot\v e\v sovec gave this asymptotic formula as well, without proof. 
If 
$$
s_n:=[q^n]\prod_{k\ge1}(1-q^{k^2})
$$ is the corresponding number for distinct squares, with the help of MAPLE we get the values
$$
(s_n)_{n=0}^{15}=(1,\,-1,\,0,\,0,\,-1,\,1,\,0,\,0,\,0,\,-1,\,1,\,0,\,0,\,1,\,-1,\,1,\,0)
$$
or
$$
(s_n)_{n=2990}^{3000}=(111,\,-112,\,61,\,46,\,-114,\,116,\,-21,\,11,\,-30,\,-17,\,37)
$$
and $\max_{n\le3000}|s_n|=319$. It is \cite[A276516]{oeis}.

\begin{prob}\tec\label{probl_on_sqcanc}
Is  $(s_n)$ unbounded? 
\end{prob}

We finish with the generalization of Example 5 to the numbers $q^{\pm,l}(n)$, $l\in\N$. These result from an almost
complete cancellation because by the pentagonal identity, 
$$
\left|q^{\pm,l}(n)\right|\le[q^n]\bigg(\sum_{n=0}^{\infty}
\left|q^{\pm,1}(n)\right|q^n\bigg)^l\le[q^n](1-q)^{-l}=\binom{n+l-1}{l-1}\;.
$$
So $q^{\pm,l}(n)=O(n^{l-1})$, $m(n)=\Theta_l(\log(1+n))$ for this counting problem and effective computation of 
$q^{\pm,l}(n)$ means computation in $\mathrm{poly}(\log n)$ steps. Besides the pentagonal identity for $l=1$, 
another nice identity occurs for $l=3$: 
$$
\prod_{n\ge1}(1-q^n)^3=\sum_{n\ge0}(-1)^{n}(2n+1)q^{n(n+1)/2}\;,
$$ 
due to C. Jacobi (G.\,H. Hardy and E.\,M. Wright \cite[Theorem 357]{hard_wrig}). Thus also $q^{\pm,3}(n)$ is a PIO 
function. For $l=2$, Example 5, we get
$$
(q^{\pm,2}(n))_{n\ge0}=(1,\,-2,\,-1,\,2,\,1,\,2,\,-2,\,0,\,-2,\,-2,\,1,\,0,\,0,\,2,\,3,\,-2,\,2,\,0,\,\ds)
$$
or
$$
(q^{\pm,2}(n))_{n=58}^{75}=(0,\,-2,\,0,\,-2,\,0,\,-2,\,2,\,0,\,-4,\,0,\,0,\,-2,\,-1,\,2,\,0,\,2,\,0,\,0)\;,
$$
\cite[A002107]{oeis}, not showing any clear pattern. J.\,W.\,L. Glaisher \cite[p. 183]{glai} writes: 
``I had no hope that these 
coefficients would follow any simple law, as in the Eulerian or Jacobian series; for, if such a law existed, it 
could not fail to have been discovered long ago by observation.'' Let us see how we advanced in 130 years.
In August 2018 the ``links'' section of the entry \cite[A002107]{oeis} (author N.\,J.\,A. Sloane) lists these references: 
table of first 10000 values by S. Manyama, G.\,E. Andrews \cite{andr87}, M. Boylan \cite{boyl}, S. Finch \cite{finc}, 
J.\,W.\,L. Glaisher \cite{glai}, J.\,T. Joichi \cite{joic}, V.\,G. Ka\v c and D.\,H. Peterson \cite{kac_pete}, 
M. Koike \cite{koik}, V. Kot\v e\v sovec \cite{kote_int_qser}, Y. Martin \cite{mart}, T. Silverman \cite{silv}, 
index to 74 sequences in \cite{mart} by M. Somos, M. Somos \cite{somo}, and article Ramanujan Theta Functions in 
E. Weisstein \cite{weis}. J.\,W.\,L. Glaisher \cite{glai} did discover and prove a kind of simple pattern for this
sequence, which we state in the elegant form given in \cite{finc} (the other references seem not relevant for computation
of $q^{\pm,2}(n)$): $q^{\pm,2}(n)=G(12n+1)$ where $G:\N\to\Z$ is a multiplicative function, which means that 
$G(ab)=G(a)G(b)$ whenever $a,b\in\Z$ are coprime numbers, defined on prime powers $p^r$, $r\in\N$, by
$$
G(p^r)=\left\{
\begin{array}{lll}
1&\ds&p\equiv7,\,11\ (\mathrm{mod}\ 12),\,r\equiv0\ (\mathrm{mod}\ 2),\\
(-1)^{r/2}&\ds&p\equiv5\ (\mathrm{mod}\ 12),\,r\equiv0\ (\mathrm{mod}\ 2),\\
r+1&\ds&p\equiv1\ (\mathrm{mod}\ 12),\,(-3)^{(p-1)/4}\equiv1\ (\mathrm{mod}\ p),\\
(-1)^r(r+1)&\ds&p\equiv1\ (\mathrm{mod}\ 12),\,(-3)^{(p-1)/4}\equiv-1\ (\mathrm{mod}\ p),\\
0&\ds&\mbox{otherwise}\;.
\end{array}
\right.
$$
This is a PIO formula for $G(n)=G(p^r)$ if $n$ is a known prime power. But in general we do not know if 
$q^{\pm,2}(n)$ is a PIO function because we do not know how to effectively factorize numbers. For example, 
$q^{\pm,2}(58)=G(697)=G(17)G(41)=0\times\ds=0$ and $q^{\pm,2}(59)=G(709)=(-1)^1(1+1)=-2$ because $709$ is a prime 
that is $1$ modulo $12$ and $(-3)^{177}=-3^{2^7}3^{2^5}3^{2^4}3\equiv-1$ modulo $709$ as  
$3^{16}\equiv495$, $3^{32}\equiv420$ and $3^{128}\equiv29$.

For repeated parts,
$$
\sum_{n\ge0}p^{\pm,2}(n)q^n:=\prod_{n\ge1}\frac{1}{(1+q^n)^2}=\prod_{n\ge1}(1+(-q)^{2n-1})^2\;,
$$
we get
$$
(p^{\pm,2}(n))_{n\ge0}=(1,\,-2,\,1,\,-2,\,4,\,-4,\,5,\,-6,\,9,\,-12,\,13,\,-16,\,21,\,-26,\,\ds)\;,
$$
\cite[A022597]{oeis}, and see, like before, that $p^{\pm,2}(n)$ is also $(-1)^n$ times the number of partitions of 
$n$ into $2$-sorted distinct odd numbers and that we have almost no cancellation.

For $l=24$ a shift of $q^{\pm,24}(n)$ gives the {\em Ramanujan tau function $\tau(n)$},
$$
\sum_{n\ge0}\tau(n)q^n:=q\prod_{n\ge1}(1-q^n)^{24}\;.
$$
So 
$$
(\tau(n))_{n\ge1}=(1,\,-24,\,252,\,-1472,\,4830,\,-6048,\,-16744,\,84480,\,-113643,\,\ds)\;,
$$
\cite[A000594]{oeis}. Combinatorially, $\tau(n)$ is the $(-1)^{\|\la\|}$-count of partitions $\la$ of $n-1$ into parts 
in $24$-sorted $\N$ (not that this would really help for deriving properties of $\tau(n)$). Is $\tau(n)$ a PIO 
function, can we compute it effectively, in $\mathrm{poly}(\log n)$ steps (as we noted above, $\tau(n)=O(n^{23})$)? 
The Wikipedia article \cite{wiki_tau} on tau function is silent about this 
fundamental aspect but the simple and unsatisfactory answer is again that we do not know. If we could 
effectively factorize numbers, we could besides decoding secret messages also compute $\tau(n)$ effectively: 
(i) $\tau(mn)=\tau(m)\tau(n)$ if $m$ and $n$ are coprime, (ii) $\tau(p^{k+2})=\tau(p)\tau(p^{k+1})-p^{11}\tau(p^k)$
for every prime number $p$ and every $k\in\N_0$ and (iii) $\tau(p)$ can be computed in $\mathrm{poly}(\log p)$ steps
for every prime $p$. The first two properties, conjectured by S. Ramanujan, were proved by L.\,J. Mordell \cite{mord} and 
the whole book \cite{edix_couv}, edited and mostly written by B. Edixhoven and J.-M. Couveignes, is devoted to
exposition of an algorithm proving (iii). 

\begin{prob}\tec\label{probl_on_pmqnl}
Are, in the current state of knowledge, the known PIO functions $q^{\pm,l}(n)$ only those for 
$l=1$ and $3$? For which $l\in\N$ can one compute $q^{\pm,l}(n)$ in $\mathrm{poly}(\log n)$ steps with the 
help of an oracle that can factorize integers efficiently?
\end{prob}

\noindent
From the extensive literature on the numbers $q^{\pm,l}(n)$ we further mention only H.\,H. Chan, S. Cooper and 
P.\,Ch. Toh \cite{chan_al1,chan_al2} (check the former for $l=26$), E. Clader, Y. Kemper and M. Wage \cite{clad_al} 
and J.-P. Serre \cite{serr}.

\bigskip\noindent
{\bf Acknowledgments.} The OEIS database \cite{oeis} was very helpful. I thank I. Pak for valuable comments
and references.

\bigskip\noindent
{\em Martin Klazar\\
Department of Applied Mathematics\\
Charles University, Faculty of Mathematics and Physics\\
Malostransk\'e n\'am\v est\'\i\ 25\\
11800 Praha\\
Czechia\\
{\tt klazar@kam.mff.cuni.cz}
}

{\small
\begin{thebibliography}{XX}

\bibitem{abra_al}
S.\,A. Abramov, M. Petkov\v sek and H. Zakraj\v sek, Convolutions of Liouvillian sequences, ArXiv:1803.08747v1, 
2018, 24 pages.

\bibitem{AKS}
M. Agrawal, N. Kayal and N. Saxena, PRIMES is in P, {\em Ann. Math.} {\bf 160} (2004), 781--793. 

\bibitem{aign}
M. Aigner, {\em A Course in Enumeration}, Springer, Berlin, 2007.

\bibitem{alfo}
J.\,L.\,R. Alfons\'\i n, {\em The Diophantine Frobenius Problem}, Oxford University Press, Oxford, 2005. 

\bibitem{amor_zann}
F. Amoroso and U. Zannier (editors), {\em Diophantine Approximation. Lectures given at the C.I.M.E. Summer 
School held in Cetraro, Italy, June 28--July 6, 2000,} Lecture Notes in Mathematics 1819, Springer, Berlin, 2003.

\bibitem{andr}
G.\,E. Andrews, {\em The Theory of Partitions}, Addison-Wesley, Reading, MA, 1976.

\bibitem{andr87}
G.\,E. Andrews, Advanced problem 6562, {\em Amer. Math. Monthly} {\bf 94} (1987), 1011.

\bibitem{andr_beck_robb}
G.\,E. Andrews, M. Beck and N. Robbins, Partitions with fixed differences between largest and 
smallest parts, {\em Proc. Amer. Math. Soc.} {\bf 143} (2015), 4283--4289.

\bibitem{andr_newm}
G.\,E. Andrews and D. Newman, Binary representations of theta functions,  {\em Integers} {\bf 18} (2018), 
$\#A34$, 8 pages.

\bibitem{ardi}
F. Ardila, {\em Algebraic and Geometric Methods in Enumerative Combinatorics}, ArXiv:1409.2562v2, 
2014, 144 pages (also Chapter 1 in \cite{bona_handb}).

\bibitem{bach_shal}
E. Bach and J. Shallit, {\em Algorithmic Number Theory, Vol. 1: Efficient Algorithms}, The MIT Press, Cambridge, 
MA, 1996.

\bibitem{bako}
V.\,P. Bakoev, Algorithmic approach to counting of certain types $m$-ary partitions, {\em Discrete Math.} {\bf 275}
(2004), 17--41.

\bibitem{barh}
Y. Bar-Hillel (editor), {\em Logic, Methodology and Philosophy of Science. Proceedings of the 1964 
International Congress}, North-Holland, Amsterdam, 1965. 

\bibitem{bech_heib_slam}
V. Becher, P.\,A. Heiber and T.\,A. Slaman, A polynomial-time algorithm for computing absolutely normal 
numbers, {\em Information and Computation} {\bf 232} (2013), 1--9.

\bibitem{bell}
E.\,T. Bell, Interpolated denumerants and Lambert series, {\em Amer. J. Math.} {\bf 65} (1943), 382--386. 

\bibitem{berg_labe_lero}
F. Bergeron, G. Labelle and P. Leroux, {\em Combinatorial Species and Tree-like Structures}, Cambridge University 
Press, Cambridge, 1998. 

\bibitem{bodi_grop_kang}
M. Bodirsky, C. Gr\"opl and M. Kang, Generating labeled planar graphs uniformly at random, 
{\em Theoretical Comp. Sci.} {\bf 379} (2007), 377--386.

\bibitem{boga_good_nguy}
T. Bogart, J. Goodrick, D. Nguyen and K. Woods, Parametric Presburger arithmetic: complexity of counting and 
quantifier elimination, ArXiv:1802.00974v1, 2018, 14 pages.

\bibitem{boga_good_wood}
T. Bogart, J. Goodrick and K. Woods, Parametric Presburger arithmetic: logic, combinatorics, and 
quasipolynomial behavior, {\em Discrete Analysis}, 2017, 34 pages. 

\bibitem{bohm_frob_ries}
J. Bohman, C.-E. Fr\"oberg and H. Riesel, Partitions in squares, {\em BIT} {\bf 19} (1979), 297--301.

\bibitem{bomb_gubl}
E. Bombieri and W. Gubler, {\em Heights in Diophantine Geometry}, Cambridge University Press, Cambridge, 2006.

\bibitem{bona_handb}
M. B\'ona (editor), {\em Handbook of Enumerative Combinatorics}, CRC Press, Boca Raton, FL, 2015. 

\bibitem{borw}
P.\,B. Borwein, On the complexity of calculating factorials, {\em J. of Algorithms} {\bf 6} (1985), 376--380.

\bibitem{bost_caru_chri_duma}
A. Bostan, X. Caruso, G. Christol and P. Dumas, Fast coefficient computation for algebraic power series in 
positive characteristic, ArXiv:1806.06543v1, 2018, 16 pages.

\bibitem{bost_chri_duma}
A. Bostan, G. Christol and P. Dumas, Fast computation of the $N$th term of an algebraic series in positive 
characteristic, in: {\em ISSAC'16 Proceedings of the ACM on International Symposium on Symbolic and Algebraic Computation
(2016)}, 119--126. 

\bibitem{bost_gaud_scho}
A. Bostan, P. Gaudry and E. Schost, Linear recurrences with polynomial coefficients and application to
integer factorization and Cartier--Manin operator, {\em SIAM J. Comput.} {\bf 36} (2007), 1777--1806.

\bibitem{bowm_al}
D. Bowman, P. Erd\"os and A. Odlyzko, Partitions of $n$ into parts which are divisors of $n$,  {\em Amer. Math. 
Monthly} {\bf 99} (1992), 276--277.

\bibitem{boyl}
M. Boylan, Exceptional congruences for the coefficients of certain eta-product newforms, {\em J. Number Theory} 
{\bf 98} (2003), 377--389.

\bibitem{brui}
N.\,G. de Bruijn, On Mahler's partition problem, {\em Proc. Kon. Ned. Akad. v. Wet. Amsterdam} {\bf 51} (1948), 
659--669.

\bibitem{brui_ono}
J.\,H. Bruiner and K. Ono, Algebraic formulas for the coefficients of half-integral weight harmonic weak Maas forms,
{\em Adv. Math.} {\bf 246} (2013), 198--219.

\bibitem{brui_ono_suth}
J.\,H. Bruiner, K. Ono and A.\,V. Sutherland, Class polynomials for nonholomorphic modular functions,
{\em J. Number Theory} {\bf 161} (2016), 204--229.

\bibitem{calk_al}
N. Calkin, J. Davis, K. James, E. Perez and C. Swannack, Computing the integer partition function, 
{\em Math. Comp.} {\bf 76} (2007), 1619--1638.

\bibitem{came}
P.\,J. Cameron, {\em Notes on Counting: An Introduction to Enumerative Combinatorics}, 2010, 217 pages  
(available from P.\,J. Cameron's homepage).

\bibitem{came_book}
P.\,J. Cameron, {\em Notes on Counting: An Introduction to Enumerative Combinatorics}, Cambridge University Press,
Cambridge, 2017. 

\bibitem{chan_al1}
H.\,H. Chan, S. Cooper and P.\,Ch. Toh, The $26$th power of Dedekind's $\eta$-function, {\em Adv. in Math.}
{\bf 207} (2006), 532--543.

\bibitem{chan_al2}
H.\,H. Chan, S. Cooper and P.\,Ch. Toh, Ramanujan's Eisenstein series and powers of Dedekind's $\eta$-function, 
{\em J. London Math. Soc.} {\bf 75} (2007), 225--242.

\bibitem{chol_koli_sill}
Y. Choliy, L.\,W. Kolitsch and A.\,V. Sills, Partition recurrences, {\em Integers} {\bf 18B} (2018), 
article A1, 15 pages. 

\bibitem{chri}
G. Christol, Globally bounded solutions of differential equations, in: \cite{annt}, 45--64. 

\bibitem{chri_chri}
A.\,D. Christopher and M.\,D. Christober, Estimates of five restricted partition functions that are quasi 
polynomials, {\em Bull. Math. Sci.} {\bf 5} (2015), 1--11.

\bibitem{chu_al} W. Chu, G. Gardin, S. Ohsuga and Y. Kambayashi (editors), {\em Proc. 7th International Conference 
on Very Large Data Bases}, Morgan Kaufmann, Cannes, 1981. 

\bibitem{cimp_nico1}
M. Cimpoea\c{s} and F. Nicolae, On the restricted partition function, ArXiv:1609.060901v1, 2016, 20 pages
(to appear in {\em The Ramanujan J.}).

\bibitem{cimp_nico2}
M. Cimpoea\c{s} and F. Nicolae, On the restricted partition function II, ArXiv:1611.00256v4, 2018, 12 pages. 

\bibitem{ciol}
A. Ciolan, Asymptotics and inequalities for partitions into squares, ArXiv:1806.00708v1, 2018, 16 pages. 

\bibitem{clad_al}
E. Clader, Y. Kemper and M. Wage, Lacunarity of certain partition-theoretic generating functions, 
{\em Trans. Amer. Math. Soc.} {\bf 137} (2009), 2959--2968.

\bibitem{cobh}
A. Cobham, The intrinsic computational difficulty of functions, in: \cite{barh}, 24--30.

\bibitem{cohe_kime_sagi}
S. Cohen, B. Kimelfeld and Y. Sagiv, Generating all maximal induced subgraphs for hereditary and connected-hereditary
graph properties, {\em J. Computer and System Sci.} {\bf 74} (2008), 1147--1159.

\bibitem{cohe_sagi}
S. Cohen and Y. Sagiv, An incremental algorithm for computing ranked full disjunction,
{\em J. Computer and System Sci.} {\bf 73} (2007), 648--668.

\bibitem{comt}
L. Comtet, {\em Advanced Combinatorics: The Art of Finite and Infinite Expansions}, D. Reidel, Dordrecht, 1974 
(first published in French in 1970).

\bibitem{dunn_robl}
A. Dunn and N. Robles, Polynomial partition asymptotics, ArXiv:1705.00384v1, 2017, 22 pages. 

\bibitem{edga}
T. Edgar, On the number of hyper $m$-ary partitions, {\em Integers} {\bf 18} (2018), article A47.

\bibitem{edix_couv}
B. Edixhoven and J.-M. Couveignes (editors), {\em Computational aspects of modular forms and Galois representations. 
How one can compute in polynomial time the values of Ramanujan's tau at a prime}, Princeton University Press, 
Princeton, NJ, 2011. 

\bibitem{edmo}
J. Edmonds, Path, trees, and flowers, {\em Canad. J. Math.} {\bf 17} (1965), 449--467. 

\bibitem{eule}
L. Euler, Evolutio producti infiniti $(1-x)(1-xx)(1-x^3)(1-x^4)(1-x^4)(1-x^5)(1-x^6)$ etc. in seriem simplicem, 
{\em Acta Academiae Scientiarum Imperialis Petropolitanae (1780)} (1783), 47--55. 

\bibitem{ever_al}
G. Everest, A. van der Poorten, I. Shparlinski and T. Ward, {\em Recurrence Sequences}, AMS, Providence, RI, 2003. 

\bibitem{ever}
J.-H. Evertse, On sums of $S$-units and linear recurrences, {\em Compositio Math.} {\bf 53} (1984), 225--244.

\bibitem{fill_al}
J.\,A. Fill, S. Janson and M.\,D. Ward, Partitions with distinct multiplicities of parts: On an ``Unsolved 
Problem'' posed by Herbert Wilf, {\em Electronic J. Combinatorics} {\bf 19} (2012), Paper \#P18, 6 pages.

\bibitem{finc}
S. Finch, Powers of Euler's $q$-series, ArXiv:math/0701251v2, 2007, 17 pages. 

\bibitem{flaj_sedg}
P. Flajolet and R. Sedgewick, {\em Analytic Combinatorics}, Cambridge University Press, Cambridge, 2009.

\bibitem{frob}
C.-E. Fr\"oberg, Accurate estimation of the number of binary partitions, {\em BIT} {\bf 17} (1977), 386--391.

\bibitem{gafn}
A. Gafni, Power partitions, {\em J. Number Theory} {\bf 163} (2016), 19--42.

\bibitem{gath_gerh}
J. von zur Gathen and J. Gerhard, {\em Modern Computer Algebra}, Cambridge University Press, Cambridge, 
2013 (3rd edition).

\bibitem{glai1}
J.\,W.\,L. Glaisher, On formulae of verification in the partition of numbers, {\em Proc. Royal Soc. London} 
{\bf 24} (1876), 250--259.

\bibitem{glai}
J.\,W.\,L. Glaisher, On the square of Euler's series, {\em Proc. London Math. Soc.} {\bf 21} (1889), 182--215.

\bibitem{gowe}
T. Gowers, J. Barrow-Green and I. Leader (editors), {\em The Princeton Companion to Mathematics}, Princeton 
University Press, Princeton, NJ, 2008.

\bibitem{gure_shel} 
Y. Gurevich and S. Shelah, Time polynomial in input or output, {\em J. Symb. Logic} {\bf 54} (1989), 1083--1088.

\bibitem{hala}
G. Hal\'asz (editor), {\em Topics in Classical Number Theory. Vol. II (Colloq. Math. Soc. J. Bolyai 34)},
North-Holland, Amsterdam, 1984. 

\bibitem{hard_rama}
G.\,H. Hardy and S. Ramanujan, Asymptotic formulae in combinatory analysis, {\em Proc. London Math. Soc.} {\bf 2} 
(1918), 75--115.

\bibitem{hard_wrig}
G.\,H. Hardy and E.\,M. Wright, {\em An Introduction to the Theory of Numbers}, Clarendon Press, Oxford, 1979 
(fifth edition, first published in 1938).

\bibitem{harv_hoev}
D. Harvey and J. van der Hoeven, Faster integer multiplication using short lattice vectors, ArXiv:1802.07932v1, 
2018, 16 pages. 

\bibitem{harv_hoev_lece}
D. Harvey, J. van der Hoeven and G. Lecerf, Even faster integer multiplication, {\em J. Complexity} {\bf 36} (2016), 
1--30.

\bibitem{henn}
F.\,C. Hennie, One-tape, off-line Turing machine computations, {\em Information and Control} {\bf 8} (1965), 553--578.

\bibitem{heub_mans}
S. Heubach and T. Mansour, {\em Combinatorics of Compositions and Words}, Chapman and Hall/CRC, Boca Raton, FL, 2009.

\bibitem{hirs00}
M.\,D. Hirschhorn, Partial fractions and four classical theorems of number theory, {\em Amer. Math. Monthly} 
{\bf 107} (2000), 260--264.

\bibitem{hirs}
M.\,D. Hirschhorn, The number of partitions of a number into distinct squares, {\em Math. Gazette} {\bf 90} 
(2006), 80--87.

\bibitem{hirs_rama}
M.\,D. Hirschhorn, My contact with Ramanujan, {\em J. Indian Math. Soc. (N. S.)}, (2013), 33--43 
(special volume to commemorate the 125th birth anniversary of Srinivasa Ramanujan).

\bibitem{hirs_sell}
M.\,D. Hirschhorn and J.\,A. Sellers, A diferent view of $m$-ary partitions, {\em Australasian J. Combinatorics}
{\bf 30} (2004), 193--196. 

\bibitem{hyun_al}
S.\,G. Hyun, S. Melczer and C. St-Pierre, A fast algorithm for solving linearly recurrent sequences, 
ArXiv:1806.03554v1, 2018, 4 pages. 

\bibitem{jaki}
R. Jakimczuk, Restricted partitions, {\em Internat. J. of Mathem. and Mathem. Sciences} {\bf 36} (2004), 1893--1896.

\bibitem{jeli_klaz}
V. Jel\'\i nek and M. Klazar, Generalizations of Khovanski\u\i's theorem on the growth of sumsets in Abelian
semigroups, {\em Adv. in Appl. Math.} {\bf 41} (2008), 115--132. 

\bibitem{joha}
F. Johansson, Efficient implementation of the Hardy--Ramanujan--Rademacher formula, {\em LMS J. Comput. Math.}
{\bf 15} (2012), 341--359.

\bibitem{joha_blog}
F. Johansson, New partition function record: $p(10^{20})$ computed,\\ 
{\tt fredrikj.net/blog/2014/new-partition-function-record/}, 2014 (viewed in August 2016). 

\bibitem{joic}
J.\,T. Joichi, Hecke--Rogers, Andrews identities; combinatorial proofs, {\em Discrete Math.} {\bf 84} 
(1990), 255--259.

\bibitem{kac_pete}
V.\,G. Ka\v c and D.\,H. Peterson, Infinite-dimensional Lie algebras, theta functions and modular forms, 
{\em Adv. in Math.} {\bf 53} (1984), 125--264. 

\bibitem{kane_rhoa}
D. Kane and R.\,C. Rhoades, Asymptotics for Wilf's partitions with distinct multiplicities, preprint, 2012, 
8 pages. 

\bibitem{kanz_sagi}
Y. Kanza and Y. Sagiv, Computing full disjunctions, in: {\em Proceedings of the 22nd ACM SIGMOND-SIGACT-SIGART Symposium 
on Principles of Database Systems}, ACM Press, San Diego, CA, 2003, 78--89.

\bibitem{kaue_zeil}
M. Kauers and D. Zeilberger, A simple re-derivation of Onsager's solution of the 2D Ising model using 
experimental mathematics, ArXiv:1805.09057v1, 2018, 10 pages.

\bibitem{kenn}
J. Kennedy (editor), {\em Interpreting G\"odel. Critical Essays}, Cambridge University Press, Cambridge, 2014.

\bibitem{klaz10}
M. Klazar, Overview of some general results in combinatorial enumeration, in: \cite{perm_patt}, 3-40.

\bibitem{knop_robb}
A. Knopfmacher and N. Robbins, Identities for the total number of parts in partitions of integers, {\em Util.
Math.} {\bf 67} (2005), 9--18. 

\bibitem{knut66}
D.\,E. Knuth, An almost linear recurrence, {\em Fibonacci Q.} {\bf 4} (1966), 117--128.

\bibitem{knut}
D.\,E. Knuth, {\em The Art of Computer Programming, Vol. 2: Seminumerical Algorithms}, Addison-Wesley, Reading, MA, 
1997.

\bibitem{knut_fasc}
D.\,E. Knuth, {\em The Art of Computer Programming, Pre-fascicle 3B, A Draft of Sections 7.2.1.4--5: Generating 
all Partitions}, preprint, 2004.

\bibitem{koik}
M. Koike, On McKay conjecture, {\em Nagoya Math. J.} {\bf 95} (1984), 85--89.

\bibitem{kote}
V. Kot\v e\v sovec, A method of finding the asymptotics of $q$-series based on the convolution of generating
functions, ArXiv:1509.08708v3, 2016, 28 pages.

\bibitem{kote_int_qser}
V. Kot\v e\v sovec, The integration of q-series, {\tt oeis.org/A258232/a258232\_2.pdf}, 2015, 3 pages
(viewed in August 2018).  

\bibitem{krat_mull}
C. Krattenthaler and T.\,W. M\"uller, Motzkin numbers and related sequences modulo powers of $2$, 
ArXiv:1608.05657v1, 2016, 28 pages. 

\bibitem{kren_wagn}
D. Krenn and S. Wagner, Compositions into powers of $b$: asymptotic enumeration and parameters, 
{\em Algorithmica} {\bf 75} (2016), 606--631.

\bibitem{perm_patt}
S. Linton, N. Ruskuc and V. Vatter (editors), {\em Permutation Patterns. St. Andrews 2007}, vol. 376 of London 
Mathematical Society Lecture Note Series, Cambridge Unversity Press, Cambridge, 2010. 

\bibitem{luth}
S.\,M. Luthra, On the average number of summands in partitions of $n$, {\em Proc. Nat. Inst. Sci. India, Part A} 
{\bf 23} (1957), 483--498.

\bibitem{mahl}
K. Mahler, On a special functional equation, {\em J. London Math. Soc.} {\bf 68} (1940), 115--123.

\bibitem{mart}
Y. Martin, Multiplicative $\eta$-quotients, {\em Trans. of the Amer. Math. Soc.} {\bf 348} (1996), 4825--4856.

\bibitem{mein}
G. Meinardus, \"Uber Partitionen mit Differenzenbedingungen, {\em Math. Zeit.} {\bf 61} (1954), 289--302.

\bibitem{mill_spen}
J.\,C.\,P. Miller and D.\,J. Spencer Brown, An algorithm for evaluation of remote terms in a linear recurrence 
sequence, {\em Comput. J.} {\bf 9} (1966), 188--190. 

\bibitem{mill}
V.\,S. Miller, Counting matrices that are squares, ArXiv:1606.09299v1, 2016, 37 pages. 

\bibitem{mord}
L.\,J. Mordell, On Mr. Ramanujan's empirical expansions of modular functions, {\em Proc. Cambridge Phil. Soc.}
{\bf 19} (1917), 117--124. 

\bibitem{murt_al}
M.\,V.\,N. Murthy, M. Brack, R.\,K. Bhaduri and J. Bartel, Semi-classical analysis of distinct square partitions,
ArXiv:1808.05146v1, 2018, 38 pages.  

\bibitem{annt}
K. Nagasaka and E. Fouvry (editors), {\em Analytic Number Theory}, Lectures Notes in Mathematics 1434, Springer, 
Berlin, 1990. 

\bibitem{naka}
B. Nakamura, Approaches for enumerating permutations with a prescribed number of occurrences of 
patterns, {\em Pure Math. Appl. } {\bf 24} (2013), 179--194.

\bibitem{naka_zeil}
B. Nakamura and D. Zeilberger, Using Noonan--Zeilberger functional equation to enumerate (in polynomial time!)
generalized Wilf classes, {\em Adv. in Appl. Math. } {\bf 50} (2013), 356--366.

\bibitem{ono}
K. Ono, The distribution of the partition function modulo $m$, {\em Ann. of Math.} {\bf 151} (2000), 293--307.

\bibitem{ono_modu}
K. Ono, {\em The Web of Modularity: Arithmetic of the Coefficients of Modular Forms and $q$-series}, AMS, 
Providence, RI, 2004. 

\bibitem{pak}
I. Pak, Complexity problems in enumerative combinatorics,  ArXiv:1803.06636v2, 2018, 31 pages. 

\bibitem{pak_yeli}
I. Pak and D. Yeliussizov, in preparation.

\bibitem{pak_yeli1}
I. Pak and D. Yeliussizov, Ehrhart polynomial of some Schlafli simplices, 2018 (talk at Joint Mathematics 
Meeting, San Diego, Ca, January 2018).

\bibitem{papa}
Ch.\,H. Papadimitriou, {\em Computational Complexity}, Addison-Wesley, Reading, MA, 1994. 

\bibitem{park_shan}
T.\,R. Parkin and D. Shanks, On the distribution of parity in the partition function, {\em Math. Comput.}
{\bf 21} (1967), 466--480.

\bibitem{pema_wils}
R. Pemantle and M.\,C. Wilson, {\em Analytic Combinatorics in Several Variables}, Cambridge University Press, 
Cambridge, 2013.

\bibitem{poly}
G. P\'olya, Kombinatorische Anzahlbestimmungen f\"ur Gruppen, Graphen and chemische Verbindungen, 
{\em Acta Math.} {\bf 68} (1937), 145--254.

\bibitem{poon}
B. Poonen, Undecidable problems: a sampler, in: \cite{kenn}, 211--241.

\bibitem{vdpo}
A. van der Poorten, Some problems of recurrent interest, {\em Macquarie Math. Reports}, Macquarie University,
Northridge, Australia, 81-0037 (1981).

\bibitem{vdpo84}
A. van der Poorten, Some problems of recurrent interest, in: \cite{hala}, 1265--1294.

\bibitem{vdpo_schl}
A. van der Poorten and H.\,P. Schlickewei, Additive relations in fields, {\em J. Austral. Math. Soc. Ser. A}
{\bf 51} (1991), 154--170.

\bibitem{pras}
V.\,V. Prasolov, {\em Polynomials}, Springer, Berlin, 2004 (translated from the Russian by Dimitry Leites).

\bibitem{reut}
Ch. Reutenauer, On a matrix representation for polynomially recursive sequences, {\em Electr. J. Combin.} {\bf 19} 
(2012), $\#$P36, 26 pages.

\bibitem{rich_knop}
B. Richmond and A. Knopfmacher, Compositions with distinct parts, {\em Aequat. Math.} {\bf 49} (1995), 86--97. 

\bibitem{rior}
J. Riordan, {\em Introduction to Combinatorial Analysis}, John Wiley, New York, 1958. 

\bibitem{robi_vign}
S. Robins and Ch. Vignat, Simple proofs and expressions for the restricted partition function and 
its polynomial part, ArXiv:1802.07310v1, 2018, 5 pages. 

\bibitem{rods_sell}
\O.\,J. R\o dseth and J.\,A. Sellers, Partitions with parts in a finite set, {\em Int. J. Number Theory} {\bf 2} 
(2006), 455--468.

\bibitem{rowl} E. Rowland, What is $\ds$ an automatic sequence?, {\em Notices AMS} {\bf 62} (2015), 274--276. 

\bibitem{rowl_yass} E. Rowland and R. Yassawi, Automatic congruences for diagonals of rational functions,
{\em J. Th\'eor. Nombres Bordeaux} {\bf 27} (2015), 245--288.

\bibitem{schm}
W.\,M. Schmidt, Linear recurrence sequences, in: \cite{amor_zann}, 171--247.

\bibitem{serr}
J.-P. Serre, Sur la lacunarite des puissances de $\eta$, {\em Glasgow Math. J.} {\bf 27} (1985), 203--221. 

\bibitem{shal}
J. Shallit, The Logical Approach to Automatic Sequences. Part 4: Enumeration and Automatic Sequences, 2016,\\
{\tt https://cs.uwaterloo.ca/\~{}Shallit/Talks/linz4a.pdf}\\
(viewed September 2018)

\bibitem{silv}
T. Silverman, Counting cliques in finite distant graphs, ArXiv:1612.08085v1, 2016, 16 pages. 

\bibitem{skol}
T. Skolem, Ein Verfahren zur Behandlung gewisser exponentialer Gleichungen, in: {\em Comptes rendus du congr\'es 
des math\'ematiciens scandinaves}, Stockholm, 1934 (1935), 163--188.

\bibitem{somo}
M. Somos, Introduction to Ramanujan theta functions,\\ 
{\tt https://somos.crg4.com/multiq.pdf} (viewed in August 2018).

\bibitem{spra}
R. Sprague, \"Uber Zerlegungen in ungleiche Quadratzahlen, {\em Math. Z.} {\bf 51} (1948), 289--290.

\bibitem{stanEC1}
R.\,P. Stanley, {\em Enumerative Combinatorics. Volume I}, Wadsworth $\&$ Brooks/Cole, Monterey, CA, 1986.

\bibitem{stol} M. Stoll, Rational and transcendental growth series for the higher Heisenberg groups, {\em Invent. 
Math.} {\bf 126} (1996), 85--109.

\bibitem{ston}
D.\,S. Stones, The many formulae for the number of Latin squares, {\em Electron. J. Combin.} {\bf 17}
(2010), $\#$A1, 46 pages. 

\bibitem{tao_str_ran}
T. Tao, {\em Structure and Randomness: pages from year one of a mathematical blog}, AMS, Providence, 
NJ, 2008.

\bibitem{tran_al}
M.\,N. Tran, M.\,V.\,N. Murty and R.\,K. Bhaduri, On the quantum density of states and partitioning  
an integer, {\em Ann. Phys.} {\bf 311} (2004), 204--219. 

\bibitem{vard}
M. Vardi, On the complexity of bounded-variable queries, in: {\em Proc. 14th Symposium on Principles of 
Database Systems}, ACM Press, San Jose, CA, 1995, 266--276.

\bibitem{vatt}
V. Vatter, Enumeration schemes for restricted permutations, {\em Combinatorics, Probability and 
Computing} {\bf 17} (2008), 137--159.

\bibitem{vaug}
R.\,C. Vaughan, Squares: additive questions and partitions, {\em Int. J. Number Theory} {\bf 11} (2015), 
1367--1409.

\bibitem{wilf} H.\,S. Wilf, What is an answer?, {\em Amer. Math. Monthly} {\bf 89} (1982), 289--292.  

\bibitem{wilf10}
H.\,S. Wilf, Some unsolved problems, 3 pages,\\ 
{\tt www.math.upenn.edu/\~{}wilf/website/UnsolvedProblems.pdf},\\ 
posted Dec. 13, 2010 (viewed in August 2018).

\bibitem{wrig}
E.\,M. Wright, Asymptotic partition formulae, III. Partitions into $k$-th powers, {\em Acta Math} {\bf 63} 
(1934), 143--191. 

\bibitem{yana} M. Yanakakis, Algorithms for acyclic database schemes, in: \cite{chu_al}, 82--94. 

\bibitem{zann}
U. Zannier, {\em Lecture Notes on Diophantine Analysis}, Scuola Normale Superiore, Pisa, 2009.

\bibitem{zeil}
D. Zeilberger, Enumerative and Algebraic Combinatorics, in: \cite{gowe}, 550--561.

\bibitem{zeil12}
D. Zeilberger, Using GENERATINGFUNCTIONOLOGY to enumerate distinct-multiplicity partitions, ArXiv:1201.493v1, 
2012, 6 pages.

\bibitem{zoo}
Complexity Zoo, {\tt http://complexityzoo.uwaterloo.ca} (zoo keeper: Scott Aaronson). 

\bibitem{oeis}
The On-line Encyclopedia of Integer Sequences, {\tt https://oeis.org} (founded in 1964 by N.\,J.\,A. 
Sloane).

\bibitem{wiki_tau}
{\tt https://en.wikipedia.org/wiki/Ramanujan\_tau\_function} (viewed in August 2018).

\bibitem{weis}
{\tt mathworld.wolfram.com/RamanujanThetaFunctions.html}\\ 
(WolframMathWorld, created by E. Weisstein).

\end{thebibliography}
}
\end{document}